\documentclass[11pt]{article}
\usepackage{amssymb,amsmath,fullpage}
\usepackage{graphicx}
\numberwithin{equation}{section}
\newtheorem {thm}{Theorem}[section]
\newtheorem {lem}[thm]{Lemma}

\newtheorem {prop}[thm]{Proposition}
\newtheorem {dfn}[thm]{Definition}

\newcommand{\la}{\langle}
\newcommand{\ra}{\rangle}

\newcommand{\proof}{\par\noindent{\em Proof.\ }}

\newcommand{\qed}{\hfill $\square$\par\smallskip}
\title{\bf\boldmath 
Elliptic extension of Gustafson's $q$-integral of type $G_2$
}
\author{
{\sc Masahiko Ito}\thanks{
Department of Mathematical Sciences, University of the Ryukyus, 
Okinawa 903-0213, Japan 
}\ \ 
and {\sc Masatoshi Noumi}\thanks{
Department of Mathematics, Kobe University, 
 Rokko, Kobe 657-8501, Japan
 }
}
\date{}

\begin{document}
\maketitle

\footnote[0]{
This work is supported by JSPS Kakenhi Grants (B)15H03626 and (C)18K03339. }
\begin{abstract}
The evaluation formula for an elliptic beta integral of type $G_2$ is proved. 
The integral is expressed by a product of Ruijsenaars' elliptic gamma functions, and 
the formula includes that of Gustafson's $q$-beta integral of type $G_2$ as a special limiting case as $p\to 0$. 
The elliptic beta integral of type $BC_1$ by van Diejen and Spiridonov
is effectively used in the proof of the evaluation formula. 
\end{abstract}
\section{Introduction}
The Askey--Wilson integral is a complex integral given by 
\begin{equation}
\label{eq:AW}
\frac{(q;q)_\infty}{2\,(2\pi\sqrt{-1})}
\int_{\mathbb{T}} \frac{ (x^2,x^{-2};q)_\infty}{\prod_{k=1}^4 (a_k x,a_k x^{-1};q)_\infty}
\frac{dx}{x}=\frac{(a_1a_2a_3a_4;q)_\infty}{\prod_{1\le i<j\le 4}(a_ia_j;q)_\infty},
\end{equation}
where $|a_k|<1$ $(k=1,\ldots,4)$ and $\mathbb{T}$ is the unit circle $\{x\in \mathbb{C}\,|\,|x|=1\}$ traversed in the positive direction. Hereafter, for a fixed $q\in \mathbb{C}^*$ satisfying $|q|<1$, we use the symbol $(u;q)_\infty=\prod_{\nu=1}^\infty(1-q^\nu u)$ and the abbreviation 
$(u_1,\ldots,u_m;q)_\infty=(u_1;q)_\infty\cdots(u_m;q)_\infty$.  
The infinite product on the right-hand side of \eqref{eq:AW} is expressed by a product of $q$-gamma functions. 
In this sense formula \eqref{eq:AW} 
can be regarded as a kind of beta integral, 
and in fact plays a fundamental role in the theory of 
Askey--Wilson 
$q$-orthogonal polynomials \cite{AW}.
This type of $q$-beta integrals has been extended to multiple $q$-beta integrals
in the framework of Macdonald theory of multivariable $q$-orthogonal polynomials associated with root systems.   
In this context the Askey--Wilson integral \eqref{eq:AW} is of type $BC_1$. 
In a series of pioneering works around 1990, Gustafson discovered various evaluation formulas for 
multiple $q$-beta integrals associated with root systems, including several remarkable identities  
which are not covered by the so-called Macdonald constant terms.
In the cases of non-simply laced root systems,
there are basically two types in Gustafson's multiple $q$-beta integrals, 
which are later called {\it type I} and {\it type II} 
in the context of \cite{vDS}. 
(See also \cite{Ito} for their explicit forms.)  
\par
In the last two decades, 
several elliptic extensions of the $q$-beta integrals have been studied, especially for those of type $BC_n$ by van Diejen and Spiridonov \cite{vDS}, Spiridonov \cite{S}, Rains \cite{Ra2010}.  
They include the elliptic extension of \eqref{eq:AW}
\begin{equation}
\label{eq:eAW}
\frac{(p;p)_\infty(q;q)_\infty}{2\,(2\pi\sqrt{-1})}
\int_{\mathbb{T}} \frac{\prod_{k=1}^6 \Gamma(a_k x,a_k x^{-1};p,q)}{\Gamma(x^2,x^{-2};p,q)}
\frac{dx}{x}=\prod_{1\le i<j\le 6}\Gamma(a_ia_j;p,q),
\end{equation}
where $|a_k|<1$ $(k=1,\ldots,6)$, under the balancing condition $a_1\cdots a_6=pq$. 
Here, for fixed $p$, $q\in \mathbb{C}^*$ satisfying $|p|<1$, $|q|<1$, 
we denote by 
$\Gamma(u;p,q)$ $(u\in\mathbb{C}^\ast)$ 
the {\em Ruijsenaars elliptic gamma function} defined by 
\begin{equation}
\label{eq:def RGamma}
\Gamma(u;p,q)=\frac{(pqu^{-1};p,q)_\infty}{(u;p,q)_\infty},\quad\mbox{where}\quad
(u;p,q)_\infty=\prod_{\mu,\nu=0}^{\infty}(1-p^\mu q^\nu u).
\end{equation}
We also use the notation $\Gamma(u_1,\ldots,u_m;p,q)=\Gamma(u_1;p,q)\cdots\Gamma(u_m;p,q)$. 
Note that $\Gamma(u;p,q)$ satisfies
\begin{equation}\label{eq:RGamma1}
\Gamma(qu;p,q)=\theta(u;p)\Gamma(u;p,q)\quad\mbox{and}\quad
\Gamma(pu;p,q)=\theta(u;q)\Gamma(u;p,q),
\end{equation}
where  $\theta(u;p)=(u,p/u;p)_\infty$ is a theta function satisfying
$
\theta(pu;p)=-\theta(u;p)/u
$, and also satisfies
\begin{equation}\label{eq:RGamma2}
\Gamma(pqu^{-1};p,q) =\frac{1}{\Gamma(u;p,q)},
\quad
\frac{1}{\Gamma(u,u^{-1};p,q)}=-u^{-1}\theta(u;p)\theta(u;q). 
\end{equation}
The Askey--Wilson integral \eqref{eq:AW} is obtained from \eqref{eq:eAW} as a special case, 
first by replacing $a_6$ with $ pa_6$ and by taking the limit $p\to 0$ and $a_5\to 0$ consecutively. 
\par
Compared with the development in the cases of classical root systems, 
the elliptic extensions of the cases of exceptional root systems are not fully studied yet. 
The aim of this paper is to prove an elliptic extension of the following 
$q$-integral formula of type $G_2$ (of type I) due to Gustafson \cite[p.\,101, Theorem 8.1]{Gu1994} and \cite{Gu1990}.
\begin{prop}[Gustafson] Suppose that $a_k \in \mathbb{C}^*\,(1\le k\le 4)$ satisfy $|a_k|<1$. Then we have 
\begin{equation}
\label{eq:Gustafson G2}
\begin{split}
&
\frac{(q;q)_\infty^2}{12\,(2\pi\sqrt{-1})^2}
\int\!\!\!\!\int_
{\mathbb{T}^2}\ 
\frac{\prod_{1\le i<j\le 3}(x_ix_j,x_i^{-1}x_j,x_ix_j^{-1},x_i^{-1}x_j^{-1};q)_\infty}
{\prod_{i=1}^3\prod_{k=1}^4(a_k x_i,a_k x_i^{-1};q)_\infty}
\frac{dx_1}{x_1}\frac{dx_2}{x_2}\\
&=
\frac{(a_1^2a_2^2a_3^2a_4^2;q)_\infty}{(a_1a_2a_3a_4;q)_\infty}
\prod_{i=1}^4\frac{(a_i;q)_\infty}{(a_i^2;q)_\infty}
\prod_{1\le i< j\le 4}\frac{1}{(a_ia_j;q)_\infty}
\prod_{1\le i< j<k\le 4}\frac{1}{(a_ia_ja_k;q)_\infty},
\end{split}
\end{equation}
where $x_3=x_1^{-1}x_2^{-1}$ and $\mathbb{T}^2$ is the 2-dimensional torus given by  
\begin{equation*}
\mathbb{T}^2=
\{(x_1,x_2)\in(\mathbb{C}^*)^2\,|\, |x_i|=1
\ (i=1,2)\}.
\end{equation*}
\end{prop}
Our main result is 
\begin{thm} 
\label{thm:eG2}
Suppose that $x_1x_2 x_3=1$ and $a_k \in \mathbb{C}^*\,(1\le k\le 5)$ satisfy $|a_k|<1$. 
Under the balancing condition $(a_1a_2a_3a_4a_5)^2=pq$, we have 
\begin{equation}
\label{eq:eG2}
\begin{split}
&
\frac{(p;p)_\infty^2(q;q)_\infty^2}{12\,(2\pi\sqrt{-1})^2}
\int\!\!\!\!\int_
{\mathbb{T}^2}\ 
\frac{\prod_{i=1}^3\prod_{k=1}^5\Gamma(a_k x_i,a_k x_i^{-1};p,q)}
{\prod_{1\le i<j\le 3}\Gamma(x_ix_j,x_i^{-1}x_j,x_ix_j^{-1},x_i^{-1}x_j^{-1};p,q)}
\frac{dx_1}{x_1}\frac{dx_2}{x_2}\\
&=
\prod_{i=1}^5\frac{\Gamma(a_i^2;p,q)}{\Gamma(a_i;p,q)}
\prod_{1\le i<j\le 5}
\Gamma(a_ia_j;p,q)
\prod_{\substack{1\le i<j\\[1pt]<k\le 5}}
\Gamma(a_ia_ja_k;p,q)
\!\prod_{\substack{1\le i<j\\[1pt]<k<l\le 5}}\!
\Gamma(a_ia_ja_ka_l;p,q).
\end{split}
\end{equation}
\end{thm}
In 2007 this formula was communicated as a conjecture by one of the authors (M. Ito) to V.~P.~Spiridonov. 
This conjecture was formulated by Spiridonov and Vartanov \cite{SV2011,SV2010} 
in the context of 
duality of {\it superconformal indices}. 
As far as we know, however, no proof of this formula has been given so far. \\
\par
\noindent
{\bf Remark 1.} Gustafson's formula \eqref{eq:Gustafson G2} is included in \eqref{eq:eG2} as a limiting case; 
first replace $a_5$ with $p^{1\over 2}a_5$, and then take the limit $p\to 0$. \\
\par
\noindent
{\bf Remark 2.} By \eqref{eq:RGamma2}, 
under the condition 
$a_1a_2a_3a_4a_5=\epsilon p^{1\over 2}q^{1\over 2}
$, where $\epsilon\in \{1,-1\}$, the right-hand side of \eqref{eq:eG2} is also expressed as  
$$
\prod_{i=1}^5\frac{\Gamma(a_i^2;p,q)}{\Gamma(a_i;p,q)\Gamma(\epsilon p^{1\over 2}q^{1\over 2}a_i;p,q)}
\prod_{1\le i<j\le 5}\frac{\Gamma(a_ia_j;p,q)}{\Gamma(\epsilon p^{1\over 2}q^{1\over 2}a_ia_j;p,q)}.
$$
This coincides with the expression in the conjecture \cite[p.\,213, (36)]{SV2010} when $\epsilon=1$. 
Moreover, by the  property 
$\Gamma(u^2;p,q)
=\Gamma(u,-u, p^{1\over 2}u,- p^{1\over 2}u, q^{1\over 2}u,- q^{1\over 2}u,p^{1\over 2}q^{1\over 2}u,- p^{1\over 2}q^{1\over 2}u;p,q)$, 
the right-hand side of \eqref{eq:eG2} is rewritten into 
\begin{equation}
\label{eq:RHS of eG2}
\prod_{i=1}^5
\Gamma(-a_i, p^{1\over 2}a_i,- p^{1\over 2}a_i, q^{1\over 2}a_i,- q^{1\over 2}a_i,- \epsilon p^{1\over 2}q^{1\over 2}a_i;p,q)
\!\!\prod_{1\le i<j\le 5}\!\!
\Gamma(a_ia_j;p,q)
\!\prod_{\substack{1\le i<j\\[1pt]<k\le 5}}\!
\Gamma(a_ia_ja_k;p,q).
\end{equation}
\par
\noindent
{\bf Remark 3.}
By the constraint $x_1x_2x_3=1$, the integrand 
\begin{equation}
\label{eq:Phi(x)1}
\Phi(x)=\frac{\prod_{i=1}^3\prod_{k=1}^5\Gamma(a_k x_i,a_k x_i^{-1};p,q)}
{\prod_{1\le i<j\le 3}\Gamma(x_ix_j,x_i^{-1}x_j,x_ix_j^{-1},x_i^{-1}x_j^{-1};p,q)}
\end{equation}
of the left-hand side of \eqref{eq:eG2} is also expressed as 
\begin{equation}
\label{eq:Phi(x)2}
\Phi(x)=\prod_{i=1}^3\frac{\prod_{k=1}^5\Gamma(a_k x_i,a_k x_i^{-1};p,q)}
{\Gamma(x_i,x_i^{-1};p,q)}
\prod_{1\le j<k\le 3}\frac{1}
{\Gamma(x_jx_k^{-1},x_kx_j^{-1};p,q)}.
\end{equation}
This expression consists of 
two parts, one depending on the short roots $\{x_i,x_i^{-1}\,|\, 1\le i\le 3\}$ 
and the other on the long roots $\{x_jx_k^{-1},x_kx_j^{-1}\,|\, 1\le j<k\le 3\}$. 
In the proof of \eqref{eq:eG2} we will use the coordinates $(z_1,z_2)$ associated with the simple roots, defined as 
\begin{equation}
\label{eq: z1=x1/x2.z2=x2}
z_1=x_1/x_2,\quad  z_2=x_2. 
\end{equation}
\par
\medskip
Theorem \ref{thm:eG2} will be proved in two steps. 
The first step is to show that both sides of \eqref{eq:eG2} satisfy a common system of $q$-difference equations,
so that we can consequently confirm both sides coincide up to a constant. 
The second step is to analyze 
asymptotic behaviors of both sides at a singularity 
in order to determine the constant. 
This method is also applicable to other elliptic beta integrals. 
In particular, we refer to \cite{INIntBC} for the $BC_n$ case including the formula \eqref{eq:eAW}, 
which might be simpler than the $G_2$ case of this paper.
\par\medskip
This paper is organized as follows. 
After defining basic terminology of the root system $G_2$ in Section \ref{section:2},
we first present in Section \ref{section:3} the explicit forms of the $q$-difference equations which the integral \eqref{eq:eG2} satisfies (Proposition \ref{prop:qDE}). 
In Section \ref{section:4} we study the analytic continuation of the integral \eqref{eq:eG2} as a meromorphic function
of the parameters $a_1,\ldots,a_4$ in a specific domain. We use this argument to show
that the integral \eqref{eq:eG2} is expressed as a product of elliptic gamma functions up to a
constant. 
In Section \ref{section:5} we explain a fundamental method, which corresponds to {\em integration by parts} in calculus, 
to deduce the $q$-difference equations for the contour integral \eqref{eq:eG2}. 
This method is formulated in terms of 
a $q$-difference coboundary operator $\nabla_{\rm sym}:\mathcal{G}_\epsilon\to \mathcal{F}$, 
where $\mathcal{G}_\epsilon$ and $\mathcal{F}$ are defined as spaces of theta functions specified by individual quasi-periodicities.   
Section \ref{section:6} is a technical part; 
we investigate in detail the source and target spaces $\mathcal{G}_\epsilon$, $\mathcal{F}$ of the operator $\nabla_{\rm sym}$. We apply this argument to proving Lemma \ref{lem:qDE1}, 
which we used to derive the $q$-difference equations in Proposition \ref{prop:qDE}. 
Section \ref{section:7} is devoted to 
asymptotic analysis of the contour integral \eqref{eq:eG2} along the singularity $a_1a_2=1$.
It is used to 
determine the explicit value of the constant, which was indefinite at the stage of Section \ref{section:4}. 
It should be noted that the elliptic beta integral \eqref{eq:eAW} of type $BC_1$ 
naturally arises in the process 
of calculation of the asymptotic behavior. 
\par
Lastly, 
we comment on our calculation of $\nabla_{\rm sym}$. 
For theta functions $\varphi\in\mathcal{G}_\epsilon$ we need to expand $\nabla_{\rm sym}\varphi\in \mathcal{F}$
as a linear combination of theta functions which belong to a particular basis of $\mathcal{F}$. 
In this paper we made use of the basis of $\mathcal{F}$ that consists of 
the Lagrange interpolation functions associated with 
the specific points ${\rm p}_{23}$, ${\rm p}_{13}$, ${\rm p}_{12}$ and ${\rm p}_{12}^*\in (\mathbb{C}^*)^2$ 
defined in Section \ref{section:6}. 
(We constructed this basis in a heuristic way. See the set of theta functions $\{F_1(z),F_2(z),F'_3(z),G(z)\}$, 
whose interpolation property is presented in the table below \eqref{eq:F'3(p12),F'3(p12*)}.)  
In the cases of $A_n$ and $BC_n$ root systems 
in \cite{IN2018} and \cite{IN2019,INSumBC,INIntBC,INSlaterBC}, respectively, 
we remark that Lagrange interpolation functions 
in a space of theta functions of particular quasi-periodicity are introduced 
by systematically specifying a set of reference points in $(\mathbb{C}^*)^n$. 
It would be an interesting problem to find 
a universal way which produces adequate interpolation bases for general root systems.
\section{Root system $G_2$}\label{section:2}
Let $\{\varepsilon_1,\varepsilon_2,\varepsilon_3\}$ be the standard basis of $\mathbb{R}^3$ 
with the inner product $(\cdot,\cdot)$ satisfying $(\varepsilon_i,\varepsilon_j)=\delta_{ij}$, and let 
$V$ be the hyperplane in $\mathbb{R}^3$ with equation $\xi_1+\xi_2+\xi_3=0$, i.e., 
$
V=\{\xi\in \mathbb{R}^3\,|\, (\xi, \varepsilon_1+\varepsilon_2+\varepsilon_3)=0\}.
$\\[-4pt]
\par
\begin{figure}[htbp]
\begin{center}
\includegraphics[height=130pt]{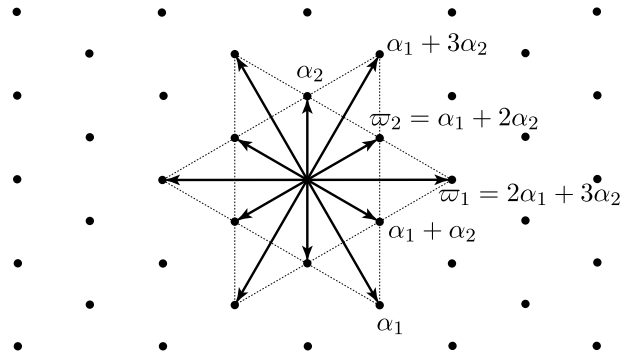}
\end{center}
\vspace{-10pt}
\caption{Root system $R$ and lattice points in $P$}
\end{figure}
\noindent 
Let $R\subset V$ be the root system of type $G_2$  given by 
$$R=\{\pm\bar\varepsilon_1,\pm\bar\varepsilon_2,\pm\bar\varepsilon_3\}
\cup\{\pm(\bar\varepsilon_1-\bar\varepsilon_2),
\pm(\bar\varepsilon_1-\bar\varepsilon_3),\pm(\bar\varepsilon_2-\bar\varepsilon_3)\},
$$
where $\bar\varepsilon_i=\varepsilon_i-(\varepsilon_1+\varepsilon_2+\varepsilon_3)/3$. 
We refer the setting of the root system of type $G_2$ to Macdonald's book \cite{Mac}.
We fix the set of simple roots $\{\alpha_1,\alpha_2\}\subset R$ given by 
$$\alpha_1=\bar\varepsilon_1-\bar\varepsilon_2=\varepsilon_1-\varepsilon_2,\qquad
\alpha_2=\bar\varepsilon_2=(-\varepsilon_1+2\varepsilon_2-\varepsilon_3)/3.$$
The set of positive roots is given by 
\begin{align*}
R^+
&=\{\bar\varepsilon_1,\bar\varepsilon_2,-\bar\varepsilon_3\}
\cup\{\bar\varepsilon_1-\bar\varepsilon_2,\bar\varepsilon_1-\bar\varepsilon_3,\bar\varepsilon_2-\bar\varepsilon_3\}\\
&=\{\alpha_2,\alpha_1+\alpha_2,\alpha_1+2\alpha_2\}\cup\{\alpha_1,\alpha_1+3\alpha_2,2\alpha_1+3\alpha_2\}.
\end{align*}
We also fix the set of fundamental weights 
$\{\varpi_1,\varpi_2\}$ by 
$(\alpha_i^\vee,\varpi_j)=\delta_{ij}$, 
where $\alpha^\vee=2\alpha/(\alpha,\alpha )$. 
This implies that 
$$\varpi_1=2\alpha_1+3\alpha_2,\qquad\varpi_2=\alpha_1+2\alpha_2.$$
Let $P$ and $Q$ be the weight lattice and root lattice defined by 
$
P=\mathbb{Z}\varpi_1+\mathbb{Z}\varpi_2
$
and 
$
Q=\mathbb{Z}\alpha_1+\mathbb{Z}\alpha_2$, respectively. 
For the root system $G_2$, the root lattice $Q$ coincides with the weight lattice $P$. 

Let $W$ be the Weyl group of type $G_2$ generated by orthogonal reflections $s_\alpha$ $(\alpha\in R)$
with respect to the hyperplane perpendicular to $\alpha\in R$,  
which are given by $s_\alpha(\xi)=\xi-(\alpha^\vee, \xi)\alpha$. 
The group $W$ is generated by the reflections $s_i=s_{\alpha_i}:V\to V$ $(i=1,2)$,  
and is isomorphic to the dihedral group of order 12. 
\[
\renewcommand{\arraystretch}{1.2}
\begin{array}{|c|cc|}
\hline
  & \alpha_1  & \alpha_2  \\
 \hline 
s_1 &  -\alpha_1 &  \alpha_1+\alpha_2 \\[3pt]
s_2 &  \alpha_1+3\alpha_2 &  -\alpha_2 \\[3pt]
\hline
\end{array}
\]
Moreover, 
$W$ is explicitly written as 
\begin{equation}
\label{eq:W}
W=\{(s_1s_2)^k,(s_1s_2)^ks_2\,|\, k=0,1,\ldots,5\},
\end{equation}
where
$(s_1s_2)^k$ coincide with the rotations around the origin through angle $k\pi/3$ on $V$ and 
$(s_1s_2)^k s_2$ coincide with the reflections $s_\alpha$ $(\alpha\in R^+)$ written as follows:
$(s_1s_2)^0s_2=s_{\alpha_2}$, 
$(s_1s_2)^1s_2=s_{\alpha_1}$, 
$(s_1s_2)^2s_2=s_{\alpha_1+\alpha_2}$, 
$(s_1s_2)^3s_2=s_{2\alpha_1+3\alpha_2}$,
$(s_1s_2)^4s_2=s_{\alpha_1+2\alpha_2}$,
$(s_1s_2)^5s_2=s_{\alpha_1+3\alpha_2}$.
We use the expression \eqref{eq:W} of $W$ later. 
The element $w_0=(s_1s_2)^3$ is the {\em longest element} of $W$. 
Note also that 
the inner product and the reflections are uniquely extended linearly to 
$V_{\mathbb{C}}=\mathbb{C}\otimes_{\mathbb{R}}V$. 

We fix the set of fundamental coweights 
$\{\omega_1,\omega_2\}$ by 
$(\omega_i,\alpha_j)=\delta_{ij}$,  so that 
$
\omega_1=\varpi_1
$, 
$\omega_2=3\varpi_2.
$
Let $P^\vee$ be the coweight lattice defined by 
$P^\vee=\mathbb{Z}\omega_1+\mathbb{Z}\omega_2$. 
For $c\in \mathbb{C}$ and  $\omega\in P^{\vee}$ 
we denote by $S_{c, \omega}$ the $c$-shift operator with respect to $\omega$ for functions $f(\zeta)$ on $V_{\mathbb{C}}$ by 
\begin{equation}
\label{eq:c-shift wrt omega}
S_{c, \omega}f(\zeta)=f(\zeta+c\,\omega).
\end{equation}
We also define action of the Weyl group $W$ on $f(\zeta)$ by 
\begin{equation}
\label{eq:w-action}
w.f(\zeta)=f(w^{-1}\zeta)\quad (w\in W).
\end{equation}

We consider the mapping from $V_{\mathbb{C}}$ to $(\mathbb{C}^*)^2$ by 
\begin{equation}
\label{eq:mapping}
\zeta\mapsto z=(e^{2\pi \sqrt{-1}(\zeta,\alpha_1)}, e^{2\pi \sqrt{-1}(\zeta,\alpha_2)}). 
\end{equation}
If we write $\zeta\in V_{\mathbb{C}}$ with the fundamental coweights by  
$\zeta=\zeta_1\omega_1+\zeta_2\omega_2$, 
then the above mapping is written as $\zeta\mapsto z=(e^{2\pi \sqrt{-1}\zeta_1},e^{2\pi \sqrt{-1}\zeta_2})$. 
For $\lambda\in P$,  we write $z^\lambda=e^{2\pi \sqrt{-1}(\zeta,\lambda)}$. In particular, 
for $\lambda=\lambda_1\alpha_1+\lambda_2\alpha_2\in Q=P$
we have the expression $z^\lambda=z_1^{\lambda_1}z_2^{\lambda_2}$, where $z_i=z^{\alpha_i}$. 
Through \eqref{eq:w-action} and \eqref{eq:mapping},
for $w\in W$ we can define $w.z^\lambda=z^{w\lambda}$, i.e., 
$$
w.z^\lambda=w.e^{2\pi \sqrt{-1}(\zeta,\lambda)}
=e^{2\pi \sqrt{-1}(w^{-1}\zeta,\lambda)}=e^{2\pi \sqrt{-1}(\zeta,w\lambda)}=z^{w\lambda},$$
and we can also define $w.f(z)$ for functions $f(z)=f(z_1,z_2)$ on $(\mathbb{C}^*)^2$ as  
$$w.f(z)=f(w.z_1,w.z_2)=f(z^{w\alpha_1},z^{w\alpha_2}),$$
so that, for instance, we have  
\begin{equation}
\label{eq:s1,s2}
s_1.f(z_1,z_2)=f(z_1^{-1},z_1z_2), \quad s_2.f(z_1,z_2)=f(z_1z_2^3,z_2^{-1}),  
\end{equation}
and 
\begin{equation}
\label{eq:w0}
w_0.f(z_1,z_2)=(s_1s_2)^3.f(z_1,z_2)=f(z_1^{-1},z_2^{-1}).
\end{equation}
We say that a function $f(z)$ is $W$-symmetric if $w.f(z)=f(z)$ for all $w\in W$. 
By chain rule for differential forms, we have 
\begin{equation}
\label{eq:dz/z}
\frac{d(s_1.z_1)}{s_1.z_1}=-\frac{dz_1}{z_1}, \ 
\frac{d(s_1.z_2)}{s_1.z_2}=\frac{dz_1}{z_1}+\frac{dz_2}{z_2}
\ \ \mbox{and}\ \ 
\frac{d(s_2.z_1)}{s_2.z_1}=\frac{dz_1}{z_1}+3\frac{dz_2}{z_2}, \ 
\frac{d(s_2.z_2)}{s_2.z_2}=-\frac{dz_2}{z_2}.
\end{equation}
\par
We fix 
$p=e^{2\pi\sqrt{-1}\sigma}$ and $q=e^{2\pi\sqrt{-1}\tau}$, where ${\rm Im}\,\sigma>0$ and 
${\rm Im}\,\tau>0$, respectively. 
If we consider a function $f(z)=f(z_1,z_2)$ on $(\mathbb{C}^*)^2$ as the function on $V_{\mathbb{C}}$ 
through \eqref{eq:mapping}, 
then the $p$-shift operators for $f(z)$ with respect to $z_i$ $(i=1,2)$ 
$$
T_{p,z_1}f(z)=f(pz_1,z_2),\qquad T_{p,z_2}f(z)=f(z_1,pz_2)
$$
are induced by the $\sigma$-shift operators $S_{\sigma, \omega_i}$ with respect to $\omega_i\in P^\vee$ $(i=1,2)$, respectively. 
The $q$-shift operators for $f(z)$ with respect to $z_i$ $(i=1,2)$ 
are also defined by  the $\tau$-shift operators with respect to $\omega_i$ $(i=1,2)$, respectively. 
\par
Using the notation 
$
x_i=e^{2\pi \sqrt{-1}(\zeta,\bar\varepsilon_i)} 
$ $(i=1,2,3)$, 
we have $x_1x_2x_3=e^{2\pi \sqrt{-1}(\zeta,\bar\varepsilon_1+\bar\varepsilon_2+\bar\varepsilon_3)}=1$ and 
the variable change $(z_1,z_2)\mapsto (x_1,x_2)$ of $(\mathbb{C}^*)^2$, where 
\begin{equation}
\label{eq:dx/x}
x_1=z_1z_2,\quad x_2=z_2
\quad%
\mbox{and}\quad
\frac{dx_1}{x_1}=\frac{dz_1}{z_1}+\frac{dz_2}{z_2},\quad
\frac{dx_2}{x_2}=\frac{dz_2}{z_2},
\end{equation}
which we saw in \eqref{eq: z1=x1/x2.z2=x2}. 
Though using the coordinates $(x_1,x_2,x_3)$ with $x_1x_2x_3=1$ instead of $(z_1,z_2)$ 
we sometimes have simple expressions for functions on $(\mathbb{C}^*)^2$ in appearance, 
like the integrands shown in \eqref{eq:Phi(x)1} or \eqref{eq:Phi(x)2} for instance, 
we use the coordinates $(z_1,z_2)$ of $(\mathbb{C}^*)^2$ associated with simple roots in the succeeding sections. 
\section{\bf\boldmath $G_2$ elliptic Gustafson integral and its $q$-difference equations}\label{section:3}
Let $\Phi(z)$ be function in $z=(z_1,z_2)\in (\mathbb{C}^*)^2$ defined by 
\begin{equation}
\label{eq:Phi}
\Phi(z)=\Phi_+(z)\Phi_+(z^{-1}),
\end{equation}
where $z^{-1}=(z_1^{-1},z_2^{-1})$ and 
$$
\Phi_+(z)=\frac{\prod_{k=1}^5\Gamma(a_k z_2,a_k z_1 z_2, a_k z_1z_2^2;p,q)}
{\Gamma(z_2,z_1z_2,z_1z_2^2,z_1,z_1z_2^3,z_1^2z_2^3;p,q)}
$$
with complex parameters $a=(a_1,\ldots,a_5)\in \mathbb{C}^*$. 
We also use the notation $\Phi(a;z)=\Phi(a_1,\ldots,a_5;z)$ instead of $\Phi(z)$ 
when we need to make the dependence on the parameters $a=(a_1,\ldots,a_5)$ explicit. 
Through \eqref{eq: z1=x1/x2.z2=x2}, $\Phi(z)$ coincides with 
\eqref{eq:Phi(x)1} or \eqref{eq:Phi(x)2}.
For the function $\Phi(z)$ we investigate the double integral 
$$
I=\int\!\!\!\!\int_{\sigma}\Phi(z)\varpi(z),
\quad\varpi(z)
=\varpi(z_1,z_2)
=\frac{1}{(2\pi\sqrt{-1})^2}\frac{dz_1}{z_1}\frac{dz_2}{z_2}
$$
over a 2-cycle $\sigma$. \eqref{eq:dx/x} implies $\varpi(x_1,x_2)=\varpi(z_1,z_2)$. 
\eqref{eq:dz/z} also implies
$
s_1.\varpi(z)=\varpi(z)
$, 
$
s_2.\varpi(z)=\varpi(z)
$,  
so that $\varpi(z)$ is $W$-symmetric. 
If the parameters satisfy the condition $|a_1|<1,\ldots,|a_5|<1$, then $\Phi(z)$ is holomorphic in the neighborhood of 
the 2-dimensional torus 
$$\mathbb{T}^2=\{z=(z_1,z_2)\in(\mathbb{C}^*)^2\,|\, |z_i|=1
\ (i=1,2)\},
$$
and hence the integral 
\begin{equation}
\label{eq:I(a)}
I(a)=I(a_1,a_2,a_3,a_4,a_5)=\int\!\!\!\!\int_{\mathbb{T}^2}\Phi(a;z)\varpi(z)
\end{equation}
defines a holomorphic function on the domain 
$
\{(a_1,\ldots,a_5)\in(\mathbb{C}^*)^5\,|\, |a_i|<1
\ (i=1,\ldots,5)\}. 
$
\par
We now formulate a system of $q$-difference equations for the integral $I(a_1,\ldots,a_5)$. 
Our goal is to establish the following proposition. 
We use the abbreviation $\theta(u_1,\ldots,u_m;p)=\theta(u_1;p)\cdots\theta(u_m;p)$.
\begin{prop}
\label{prop:qDE}
Suppose that $|p|<|q|$. Under the balancing condition $(a_1a_2a_3a_4a_5)^2=pq$, 
the integral $I(a)$ satisfies the system of $q$-difference equations 
\begin{equation}
\label{eq:qDE}
\begin{split}
I(a_1,\ldots,qa_k,\ldots,a_4,q^{-1}a_5)
&= I(a_1,\ldots,a_4,a_5)\ \frac{\theta(a_k^2,qa_k^2,q^{-1}a_5,a_1a_2a_3a_4;p)}
{\theta(q^{-2}a_5^2,q^{-1}a_5^2,a_k,q^{-1}a_1a_2a_3a_4a_5/a_k;p)}\\
&\quad\times\prod_{\substack{1\le i\le 4\\i\ne k}}\frac{\theta(a_ka_i;p)}{\theta(q^{-1}a_5a_i;p)}
\prod_{\substack{1\le i<j\le 4\\i,j\in\{1,2,3,4\}\backslash\{k\}}}\frac{\theta(a_ka_ia_j;p)}{\theta(q^{-1}a_5a_ia_j;p)}
\end{split}
\end{equation}
for $k=1,\ldots,4$, provided that $|a_1|<1,\ldots,|a_4|<1$ and $|a_5|<|q|$.
\end{prop}
Note that the condition $|a_5|<|q|$ is equivalent to $|a_1\cdots a_4|>|p|^{1\over 2}/|q|^{1\over 2}$ under the balancing condition. 
We need to assume that $p$ satisfies $|p|<|q|$ to guarantee that the above equations hold in a nonempty region. 
\par
For $u,v\in \mathbb{C}^*$ let $e(u,v;p)$ the function defined by 
\begin{equation}
\label{eq:e(u,v)}
e(u,v;p)=u^{-1}\theta(uv,u/v;p),
\end{equation}
which satisfies 
$
e(u,v;p)=-e(v,u;p)$, 
$e(pu,v;p)=e(u,v;p)(pu^2)^{-1} 
$ 
and the three-term relation
\begin{equation}
\label{eq:R/W-relation}
e(u,v;p)e(w,x;p)-e(u,w;p)e(v,x;p)+e(u,x;p)e(v,w;p)=0.
\end{equation}
We use the notation 
\begin{equation}
\label{eq:la phi ra}
\la \varphi(z)\ra=\int\!\!\!\!\int_{\mathbb{T}^2} \varphi(z)\Phi(z)\varpi(z)
\end{equation}
for any meromorphic function $\varphi(z)$ on $(\mathbb{C}^*)^2$ such that $\varphi(z)\Phi(z)$ is holomorphic in a neighborhood of $\mathbb{T}^2$. 
Since 
$\Phi(z)$ satisfies that
\begin{equation}
T_{q,a_k}\Phi(z)=a_k^3F_k(z)\Phi(z)\quad (k=1,\ldots,5),
\label{eq:TakPhi(z)}
\end{equation}
where
\begin{equation}
F_k(z)=e(a_k,z_2;p)e(a_k,z_1z_2;p)e(a_k,z_1z_2^2;p)\quad (k=1,\ldots,5),
\label{eq:F_k(z)}
\end{equation}
the integral 
$
I(a)=\la 1\ra
$
satisfies  
$$
T_{q,a_k}I(a)=a_k^3\la F_k(z)\ra \quad (k=1,\ldots,5). 
$$
\begin{lem}
\label{lem:qDE1}
Under the condition $|p|<|q|$, $(a_1a_2a_3a_4a_5)^2q=p$ and $|a_i|<1$ $(i=1,\ldots,5)$, we have  
\begin{equation}
\label{eq:F1/F2}
\la F_1(z) \ra
=\la F_2(z) \ra
\frac{a_2^3}{a_1^3}
\frac{\theta(a_1^2,qa_1^2,a_2,a_1a_3a_4a_5;p)}{\theta(a_2^2,qa_2^2,a_1,a_2a_3a_4a_5;p)}
\prod_{i=3}^5\frac{\theta(a_1a_i;p)}{\theta(a_2a_i;p)}
\prod_{3\le i<j\le 5}\frac{\theta(a_1a_ia_j;p)}{\theta(a_2a_ia_j;p)}.
\end{equation}
In other words, 
\begin{align*}
&I(qa_1,a_2,a_3,a_4,a_5)\\[-2pt]
&=I(a_1,qa_2,a_3,a_4,a_5)
\frac{\theta(a_1^2,qa_1^2,a_2,a_1a_3a_4a_5;p)}
{\theta(a_2^2,qa_2^2,a_1,a_2a_3a_4a_5;p)}
\prod_{i=3}^5\frac{\theta(a_1a_i;p)}{\theta(a_2a_i;p)}
\prod_{3\le i<j\le 5}\frac{\theta(a_1a_ia_j;p)}{\theta(a_2a_ia_j;p)}.
\end{align*}
\end{lem}
{\it Proof.} The proof of this lemma will be given later in Section \ref{section:6}. \qed
\par
\medskip
In general,  we have the following.
\begin{lem}
\label{lem:qDE2}
Under the condition $|p|<|q|$, $(a_1a_2a_3a_4a_5)^2q=p$ and $|a_i|<1$ $(i=1,\ldots,5)$, the integral $I(a)$ satisfies 
the two-term relations 
\begin{align*}
&I(a_1,\ldots,qa_k,\ldots,a_4,a_5)\\
&= I(a_1,\ldots,a_4,qa_5)
\frac{\theta(a_k^2,qa_k^2,a_5,a_1a_2a_3a_4;p)}
{\theta(a_5^2,qa_5^2,a_k,a_1a_2a_3a_4a_5/a_k;p)}
\prod_{\substack{1\le i\le 4\\i\ne k}}\frac{\theta(a_ka_i;p)}{\theta(a_5a_i;p)}
\prod_{\substack{1\le i<j\le 4\\i,j\in\{1,2,3,4\}\backslash\{k\}}}\frac{\theta(a_ka_ia_j;p)}{\theta(a_5a_ia_j;p)}
\end{align*}
for $k=1,\ldots,4$.
\end{lem}
\par
Proposition \ref{prop:qDE} is obtained from Lemma \ref{lem:qDE2} replacing $a_5$ by $q^{-1}a_5$.
\par
\medskip
We now suppose that $(a_1\cdots a_5)^2=pq$ and regard 
$a_5=\epsilon p^{1\over 2}q^{1\over 2}/a_1a_2a_3a_4$, where $\epsilon\in \{1,-1\}$, 
as a function of $(a_1,\ldots,a_4)$. 
Then the integral $I(a_1,\ldots,a_5)$ is defined on the nonempty open subset 
\begin{equation}
\label{eq:U0}
U_0=\{(a_1,\ldots,a_4)\in (\mathbb{C}^*)^4\,|\, |a_1|<1,\ldots, |a_4|<1\mbox{ and } |a_1\cdots a_4|>|p|^{1\over 2}|q|^{1\over 2}\}
\end{equation}
of $(\mathbb{C}^*)^4$. 
The $q$-difference equations \eqref{eq:qDE} for $I(a_1,\ldots,a_5)$ are defined on 
\begin{equation}
\label{eq:V0}
V_0=\{(a_1,\ldots,a_4)\in (\mathbb{C}^*)^4\,|\, |a_1|<1,\ldots, |a_4|<1\mbox{ and }  |a_1\cdots a_4|>|p|^{1\over 2}/|q|^{1\over 2}\},
\end{equation}
which is a nonempty open subset of $U_0$, if $|p|<|q|$.  
\section{Analytic continuation} \label{section:4}
The integral $I(a_1,\ldots,a_5)$, regarded as a holomorphic function in $(a_1,\ldots,a_4)\in U_0$, can be  continued to a meromorphic function on $(\mathbb{C}^*)^4$.
We prove this fact by means the $q$-difference 
equations \eqref{eq:qDE}.  
\par
In view of Proposition \ref{prop:qDE} 
we consider the meromorphic function
\begin{align}
\label{eq:defJ}
&J(a_1,\ldots,a_5)\nonumber
\\&\quad
=\prod_{i=1}^5\frac{\Gamma(a_i^2;p,q)}{\Gamma(a_i;p,q)}
\prod_{1\le i<j\le 5}
\Gamma(a_ia_j;p,q)
\prod_{\substack{1\le i<j\\[1pt]<k\le 5}}
\Gamma(a_ia_ja_k;p,q)
\!\prod_{\substack{1\le i<j\\[1pt]<k<l\le 5}}\!
\Gamma(a_ia_ja_ka_l;p,q),
\end{align}
which is also written as \eqref{eq:RHS of eG2} if $(a_1\ldots a_5)^2=pq$, as is mentioned in the introduction. 
Then it turns out that $J(a_1,\ldots,a_5)$ 
satisfies the same $q$-difference 
equations as \eqref{eq:qDE}.  
In fact, from \eqref{eq:RGamma1} one has 
\begin{equation*}
\begin{split}
J(a_1,\ldots,qa_k,\ldots,a_4,q^{-1}a_5)
&= J(a_1,\ldots,a_4,a_5)\ \frac{\theta(a_k^2,qa_k^2,q^{-1}a_5,a_1a_2a_3a_4;p)}
{\theta(q^{-2}a_5^2,q^{-1}a_5^2,a_k,q^{-1}a_1a_2a_3a_4a_5/a_k;p)}\\
&\quad\times\prod_{\substack{1\le i\le 4\\i\ne k}}\frac{\theta(a_ka_i;p)}{\theta(q^{-1}a_5a_i;p)}
\prod_{\substack{1\le i<j\le 4\\i,j\in\{1,2,3,4\}\backslash\{k\}}}\frac{\theta(a_ka_ia_j;p)}{\theta(q^{-1}a_5a_ia_j;p)}
\end{split}
\end{equation*}
for $k=1,\ldots,4$.
In the following we regard $J(a_1,\ldots,a_5)$ as 
a meromorphic function in $(a_1,\ldots,a_4)$ through 
$a_5=\epsilon p^{1\over 2}q^{1\over 2}/a_1a_2a_3a_4$, where $\epsilon\in \{1,-1\}$, as before.  
Noting that the integral 
$I(a_1,\ldots,a_5)$ is a holomorphic function on 
$U_0$, we consider the meromorphic function 
\begin{equation*}
\begin{split}
&f(a_1,\ldots,a_5)
=\frac{I(a_1,\ldots,a_5)}{J(a_1,\ldots,a_5)}
\\&\quad
=I(a_1,\ldots,a_5)
\prod_{i=1}^5\frac{1}
{\Gamma(-a_i,p^{1\over 2} a_i,-p^{1\over 2} a_i,q^{1\over 2} a_i,-q^{1\over 2} a_i,
-\epsilon p^{1\over 2}q^{1\over 2} a_i
;p,q)}\\
&\hspace{80pt}\times
\prod_{1\le i<j\le 5}
\frac{1}
{\Gamma(a_ia_j;p,q)}
\prod_{1\le i<j<k\le 5}
\frac{1}
{\Gamma(a_ia_ja_k;p,q)
}
\\[4pt]
&\quad=
I(a_1,\ldots,a_5)
\prod_{i=1}^5\frac
{(-a_i,p^{1\over 2} a_i,-p^{1\over 2} a_i,q^{1\over 2} a_i,-q^{1\over 2} a_i,
-\epsilon p^{1\over 2}q^{1\over 2} a_i;p,q)_\infty}
{(-pq/a_i,p^{1\over 2}q/a_i,-p^{1\over 2}q/a_i,
pq^{1\over 2}/a_i,-pq^{1\over 2}/a_i,-\epsilon p^{1\over 2}q^{1\over 2}/a_i
;p,q)_\infty}
\\
&\hspace{80pt}\times
\prod_{1\le i<j\le 5}
\frac
{(a_ia_j;p,q)_\infty}
{(pq/a_ia_j;p,q)_\infty}
\prod_{1\le i< j<k\le 5}
\frac
{(a_ia_ja_k;p,q)_\infty}
{(pq/a_ia_ja_k;p,q)_\infty}
\end{split}
\end{equation*}
on $U_0$.  
This ratio $f(a_1,\ldots,a_5)$ 
has poles in $U_0$ possibly along the divisors
\begin{align}
\label{eq:poles}
-a_i&=p^{\mu+1} q^{\nu+1},\  
\pm p^{1\over 2} a_i=p^{\mu+1} q^{\nu+1},\  
\pm q^{1\over 2} a_i=p^{\mu+1} q^{\nu+1},\ 
-\epsilon
p^{1\over 2}q^{1\over 2} a_i=p^{\mu+1} q^{\nu+1}\quad 
(1\le i\le 5), \nonumber\\[2pt]
a_ia_j&=p^{\mu+1} q^{\nu+1}\quad
(1\le i<j\le 5), \quad 
a_ia_ja_k=p^{\mu+1} q^{\nu+1}\quad 
(1\le i<j<k\le 5),  
\end{align}
where $\mu,\nu=0,1,2,\ldots$. 
Also, 
$f(a_1,\ldots,a_5)$ 
is $q$-periodic with respect to $(a_1,\ldots,a_4)\in V_0$ 
in the sense that 
\begin{equation*}
f(a_1,\ldots,a_5)=f(a_1,\ldots,qa_k,\ldots,q^{-1}a_5)
\end{equation*}
for $k=1,\ldots,4$.  
\begin{lem}\label{lem:W0}
Suppose that $|p|< |q|^{9}$.
Then, under the condition $(a_1a_2a_3a_4a_5)^2=pq$, there exists an open subset 
$W_0\subset (\mathbb{C}^\ast)^4$
of the form
\begin{equation}\label{eq:W0def}
W_0=\{(a_1,\ldots,a_4)\in(\mathbb{C}^\ast)^4\,|\,
s<|a_m|<1\ (1\le m\le 4)\}
\quad(0<s<|q|)
\end{equation}
such that $W_0\subset V_0$ and that 
$f(a_1,\ldots,a_5)$ is holomorphic on $W_0$.
\end{lem}
\noindent
{\em Proof.}\ 
Under the assumption $|p|< |q|^{9}$, i.e., $|p/q|^{\frac{1}{8}}< |q|$, 
one can choose positive number $s$ such that 
\begin{equation}
\label{eq:s}
|p/q|^{\frac{1}{8}}\le s<|q|.  
\end{equation}
From the condition $(a_1\cdots a_5)^2=pq$ and \eqref{eq:W0def} we have 
\begin{equation}
\label{eq:|a5|>}
|a_5|=|p^{1\over 2}q^{1\over 2}/a_1\cdots a_4|> |p^{1\over 2}q^{1\over 2}|.
\end{equation}
We first confirm that $W_0\subset V_0$. 
Suppose that $(a_1,\ldots,a_4)\in W_0$.  Then, from \eqref{eq:W0def} and \eqref{eq:s} we have  
$|a_1\cdots a_4|>s^4\ge |p/q|^{\frac{1}{2}}$, 
which means that $W_0\subset V_0$.  
We next show that $f(a_1,\ldots,a_5)$ is holomorphic in $W_0$.  
For this purpose,  
from \eqref{eq:poles} we verify the following when $(a_1,\ldots,a_4)\in W_0$: 
\begin{equation}
\label{eq:|poles|>|pq| 1}
|a_i|>|pq|,\quad |p^{1\over 2}a_i|>|pq|,\quad |q^{1\over 2}a_i|>|pq|,\quad |p^{1\over 2}q^{1\over 2}a_i|>|pq|
\quad(1\le i\le 5)
\end{equation}
and
\begin{align}
\label{eq:|poles|>|pq| 2}
|a_ia_j|>|pq|
\quad(1\le i<j\le 5),\quad
|a_ia_ja_k|>|pq|
\quad(1\le i<j<k\le 5).  
\end{align}
For \eqref{eq:|poles|>|pq| 1}, since $
|pq|^{1\over 2}<|p|^{1\over 2}<|q|^{1\over 2}<1
$,  it suffices to show 
$|p^{1\over 2}q^{1\over 2}a_i|>|pq|$ $(1\le i\le 5)$, 
which is confirmed as follows. Using \eqref{eq:s} and \eqref{eq:|a5|>} we obtain
\begin{equation*}
|p^{1\over 2}q^{1\over 2}a_i|>|p^{1\over 2}q^{1\over 2}|s
\ge |p^{1\over 2}q^{1\over 2}||p/q|^{\frac{1}{8}}=|p|^{\frac{5}{8}}|q|^{\frac{3}{8}}>|pq|
\quad\mbox{for}\quad 1\le i\le 4
\end{equation*}
and
\begin{equation*}
|p^{1\over 2}q^{1\over 2}a_5|>|p^{1\over 2}q^{1\over 2}||p^{1\over 2}q^{1\over 2}|=|pq|.
\end{equation*}
On the other hand, for \eqref{eq:|poles|>|pq| 2}, 
since $|a_ia_j|>|a_ia_ja_k|$,  it suffices to show 
$|a_ia_ja_k|>|pq|$ $(1\le i<j<k\le 5)$, which is confirmed as follows. 
Using \eqref{eq:W0def}, \eqref{eq:s} and  \eqref{eq:|a5|>} we obtain
\begin{equation*}
|a_ia_ja_k|>s^3\ge |p/q|^{3\over 8}=|p|^{\frac{3}{8}}|q|^{-\frac{3}{8}}>|pq|
\quad\mbox{for}\quad 1\le i<j<k\le 4
\end{equation*}
and 
\begin{equation*}
|a_ia_ja_5|>s^2|pq|^{1\over 2}\ge |p/q|^{1\over 4}|pq|^{1\over 2}
=|p|^{\frac{3}{4}}|q|^{\frac{1}{4}}
>|pq|
\quad\mbox{for}\quad 1\le i<j\le 4.
\end{equation*}
This completes the proof. \qed
\begin{thm}\label{thm:I=cJ}
Suppose that $|p|<|q|^{9}$. 
Under the condition $(a_1a_2a_3a_4 a_5)^2=pq$, 
the integral 
$I(a_1,\ldots,a_5)$, regarded as a holomorphic 
function in $(a_1,\ldots,a_4)\in U_0$, 
is expressed as
\begin{equation*}
I(a_1,\ldots,a_5)=b\,J(a_1,\ldots,a_5)
\end{equation*}
for some constant $b\in\mathbb{C}$  
independent of $a_1,\ldots,a_5$.  
In particular, $I(a_1,\ldots,a_5)$
is continued to 
a meromorphic function on $(\mathbb{C}^\ast)^4$. 
\end{thm}
\proof
\ By Lemma \ref{lem:W0}, there exists an open subset 
$W_0\subset(\mathbb{C}^\ast)^4$ of the form 
\eqref{eq:W0def} where $f(a_1,\ldots,a_5)$ 
$=I(a_1,\ldots,a_5)/J(a_1,\ldots, a_5)$ 
is holomorphic and satisfies 
the $q$-difference equations 
\begin{equation}\label{eq:qDEf}
f(a_1,\ldots,a_5)=f(a_1,\ldots,q a_k,\ldots,q^{-1}a_5)
\qquad(k=1,\ldots,4), 
\end{equation}
for $(a_1,\ldots,a_4)\in W_0$.  
Note that $W_0$ 
is the product of $4$ copies of an annulus 
in which the ratio of the two radii is given by 
$s<|q|$.  Hence, by the $q$-difference equations 
\eqref{eq:qDEf}, the holomorphic function 
$f(a_1,\ldots,a_5)$ on $W_0$ is continued to a 
{\em holomorphic} function on the whole 
$(\mathbb{C}^\ast)^4$.  
It must be a constant, however, since the continued 
function $f(a_1,\ldots,a_5)$ is $q$-periodic with 
respect to the variables $a_1,\ldots,a_4$.  
If we denote this constant by $b$, we 
have $I(a_1,\ldots,a_5)=b\,J(a_1,\ldots,a_5)$ 
as a holomorphic function on $U_0$, 
and hence $I(a_1,\ldots,a_5)$ is continued to a meromorphic function on $(\mathbb{C}^\ast)^4$. 
\qed 
\par
We compute the constant $b$ in Section \ref{section:7}, 
and we eventually see that $b=12/(p;p)_\infty^2(q;q)_\infty^2$. 
Once this constant has been determined, we see 
that the statement above is valid for $|p|<1$ 
without any particular restriction. 
\section{Coboundary operator $\nabla_{\rm sym}$}\label{section:5}
In this section we explain a fundamental method for 
deriving $q$-difference equations of the contour integrals 
\eqref{eq:la phi ra} based on an operator $\nabla_{\rm sym}$. 
This method corresponds to {\em integration by parts} in calculus, 
and will be used in the succeeding section for the proof of 
Lemma \ref{lem:qDE1} presented in Section \ref{section:3}.
%
%
\par
From the definition \eqref{eq:Phi} of $\Phi(z)$ we have 
\begin{equation}
\frac{T_{q,z_1}\Phi(z)}{\Phi(z)}=-\frac{f^+(z)}{T_{q,z_1}f^-(z)}, 
\label{eq:TzPhi/Phi}
\end{equation}
where 
\begin{align}
f^+(z)
&=(z_1^2z_2^3)^{-\frac{1}{2}}\frac{\prod_{k=1}^5\theta(a_kz_1z_2,a_kz_1z_2^2;p)}{\theta(z_1z_2,z_1z_2^2,z_1,z_1z_2^3,z_1^2z_2^3;p)}
=z_1^{-1}z_2^{-\frac{3}{2}}\frac{\prod_{k=1}^5\theta(a_kz_1z_2,a_kz_1z_2^2;p)}{\theta(z_1z_2,z_1z_2^2,z_1,z_1z_2^3,z_1^2z_2^3;p)},
\label{eq:f+}\\
f^-(z)&=f^{+}(z^{-1})
=z_1z_2^{\frac{3}{2}}\frac{\prod_{k=1}^5\theta(a_k z_1^{-1} z_2^{-1}, a_k z_1^{-1}z_2^{-2};p)}{\theta(z_1^{-1}z_2^{-1},z_1^{-1}z_2^{-2},z_1^{-1},z_1^{-1}z_2^{-3},z_1^{-2}z_2^{-3};p)}.
\label{eq:f-}
\end{align}
From \eqref{eq:s1,s2} and \eqref{eq:w0}, 
we can immediately confirm that
\begin{equation}
\label{eq:s_2f+=f+}
s_2.f^+(z)=f^+(z),\quad(s_1s_2)^3.f^+(z)= f^-(z). 
\end{equation}
We remark that $f^+(z)$ and $f^-(z)$ have the quasi-periodicity 
\begin{equation}
\label{eq:QP-f+}
\begin{split}
T_{p,z_1}f^+(z)&
=f^+(z)\,p(a_1\cdots a_5 )^{-2}(p^{1\over 2}z_1)^{-2}z_2^{-3},\\
T_{p,z_2}f^+(z)&
=f^+(z)\,p^{3\over 2}(a_1\cdots a_5 )^{-3}z_1^{-3}(p^{1\over 2}z_2)^{-2},\\
\end{split}
\end{equation}
\begin{equation}
\label{eq:QP-f-}
\begin{split}
T_{p,z_1}f^-(z)
&=f^-(z)\,p^{-1}(a_1\cdots a_5 )^{2}(p^{1\over 2}z_1)^{-2}z_2^{-3},\\
T_{p,z_2}f^-(z)&=f^-(z)\,p^{-{3\over 2}}(a_1\cdots a_5 )^{3}z_1^{-3}(p^{1\over 2}z_2)^{-2},
\end{split}
\end{equation}
with respect to the $p$-shifts, respectively.
\par
We denote by $\mathcal{M}((\mathbb{C}^*)^2)$ 
the $\mathbb{C}$-vector space of meromorphic functions on $(\mathbb{C}^*)^2$, 
and by $\mathcal{O}((\mathbb{C}^*)^2)$
the $\mathbb{C}$-vector space of holomorphic functions on $(\mathbb{C}^*)^2$. 
For each function $\varphi(z)\in z_{2}^{1\over 2}\mathcal{M}((\mathbb{C}^*)^2)=\{z_{2}^{1\over 2}f(z)\,|\, f(z)\in \mathcal{M}((\mathbb{C}^*)^2)\}$
we define the function $\nabla\varphi(z)$ by 
\begin{equation}
\label{eq:nabla}
(\nabla\varphi)(z)=f^+(z)T_{q,z_1}^{1\over 2}\varphi(z)+f^-(z)T_{q,z_1}^{-{1\over 2}}\varphi(z)
\in \mathcal{M}((\mathbb{C}^*)^2),
\end{equation}
and $\nabla_{\rm sym}\varphi(z)$ by the symmetrization of  $\nabla\varphi(z)$:
\begin{equation}
\label{eq:nabla_sym}
(\nabla_{\rm sym}\varphi)(z)=\sum_{w\in W}
w.(\nabla\varphi(z))
=\sum_{w\in W}
w.\Big(f^+(z)T_{q,z_1}^{1\over 2}\varphi(z)+f^-(z)T_{q,z_1}^{-{1\over 2}}\varphi(z)\Big)
\in \mathcal{M}((\mathbb{C}^*)^2)^W, 
\end{equation}
where $\mathcal{M}((\mathbb{C}^*)^2)^W$ denotes 
the $\mathbb{C}$-vector space of $W$-invariant meromorphic functions on $(\mathbb{C}^*)^2$.
\begin{lem} 
\label{lem:nabla=0}
Suppose that $|a_k|<1$ $(k=1,\ldots 5)$. 
For any $\varphi(z)\in z_{2}^{1\over 2}\mathcal{O}((\mathbb{C}^*)^2)
$, we have 
$
\la\nabla\varphi(z)\ra=0, 
$
and hence  
$
\la\nabla_{\rm sym}\varphi(z)\ra=0. 
$
\end{lem}
{\it Proof.} 
From \eqref{eq:Phi} and \eqref{eq:TzPhi/Phi}, we have 
\begin{equation}
\label{eq:Tz(Phi f-)}
T_{q,z_1}\Big(\Phi(z)f^-(z)T_{q,z_1}^{-{1\over 2}}\varphi(z)\Big)
=-\Phi(z)f^+(z)T_{q,z_1}^{{1\over 2}}\varphi(z).
\end{equation}
From \eqref{eq:f-}
we have 
\begin{equation}
\label{eq:Phi f-}
\begin{split}
&\Phi(z)f^-(z)T_{q,z_1}^{-{1\over 2}}\varphi(z)\\
&\quad=
\Phi_+(z)\frac{\prod_{k=1}^5\Gamma(a_k z_2^{-1},qa_k z_1^{-1} z_2^{-1}, qa_k z_1^{-1}z_2^{-2};p,q)}
{\Gamma(z_2^{-1},qz_1^{-1}z_2^{-1},qz_1^{-1}z_2^{-2},qz_1^{-1},qz_1^{-1}z_2^{-3},qz_1^{-2}z_2^{-3};p,q)}
z_1z_2^{3\over 2}
T_{q,z_1}^{-{1\over 2}}
\varphi(z).
\end{split}
\end{equation}
Since $z_1z_2^{3\over 2}
T_{q,z_1}^{-{1\over 2}}
\varphi(z)\in \mathcal{O}((\mathbb{C}^*)^2)$ if $\varphi(z)\in z_{2}^{1\over 2}\mathcal{O}((\mathbb{C}^*)^2)$, 
when $|a_k|<1$ $(k=1,\ldots 5)$ and $z_2$ is fixed as $|z_2|=1$, 
the right-hand side of \eqref{eq:Phi f-} as a function of $z_1$ has no poles in the annulus $|q|\le |z_1|\le1$. 
Hence, by Cauchy's integral theorem we have 
\begin{align}
&\int\!\!\!\!\int_
{\mathbb{T}^2}T_{q,z_1}\Big(\Phi(z)f^-(z)T_{q,z_1}^{-{1\over 2}}\varphi(z)\Big)\varpi(z)
\nonumber\\
&\quad=\frac{1}{(2\pi\sqrt{-1})^2}
\int_{|z_2|=1}\bigg(\int_
{|z_1|=|q|}\Phi(z)f^-(z)T_{q,z_1}^{-{1\over 2}}\varphi(z)\frac{dz_1}{z_1}\bigg)\frac{dz_2}{z_2}
\nonumber\\
&\quad=
\int\!\!\!\!\int_
{\mathbb{T}^2}
\Phi(z)f^-(z)T_{q,z_1}^{-{1\over 2}}\varphi(z)\varpi(z).
\label{eq:Cauchy}
\end{align}
Combining \eqref{eq:Tz(Phi f-)} and  \eqref{eq:Cauchy}, we have 
$$\la\nabla\varphi(z)\ra=
\int\!\!\!\!\int_
{\mathbb{T}^2}
\Phi(z)\Big(f^+(z)T_{q,z_1}^{{1\over 2}}\varphi(z)+f^-(z)T_{q,z_1}^{-{1\over 2}}\varphi(z)\Big)\varpi(z)
=0.
$$
Since $\Phi(z)$ and $\varpi(z)$ are $W$-symmetric, we therefore obtain  
\begin{align*}
\hspace{60pt}
\la\nabla_{\rm sym}\varphi(z)\ra
&=
\int\!\!\!\!\int_
{\mathbb{T}^2}\Phi(z)\nabla_{\rm sym}\varphi(z)\varpi(z)\\
&=
\int\!\!\!\!\int_
{\mathbb{T}^2}\Phi(z)\bigg(\sum_{w\in W}
w.\Big(f^+(z)T_{q,z_1}^{1\over 2}\varphi(z)+f^-(z)T_{q,z_1}^{-{1\over 2}}\varphi(z)\Big)
\bigg)\varpi(z)\\
&=
\sum_{w\in W}\int\!\!\!\!\int_
{\mathbb{T}^2}w.\bigg(\Phi(z)
\Big(f^+(z)T_{q,z_1}^{1\over 2}\varphi(z)+f^-(z)T_{q,z_1}^{-{1\over 2}}\varphi(z)\Big)
\varpi(z)\bigg)
\\
&=
\sum_{w\in W}\int\!\!\!\!\int_
{w^{-1}.\mathbb{T}^2}\Phi(z)
\Big(f^+(z)T_{q,z_1}^{1\over 2}\varphi(z)+f^-(z)T_{q,z_1}^{-{1\over 2}}\varphi(z)\Big)\varpi(z)\\
&=
\sum_{w\in W}\int\!\!\!\!\int_
{\mathbb{T}^2}\Phi(z)
\Big(f^+(z)T_{q,z_1}^{1\over 2}\varphi(z)+f^-(z)T_{q,z_1}^{-{1\over 2}}\varphi(z)\Big)\varpi(z)\\
&=0. \hspace{320pt} 
\square
\end{align*}
\begin{dfn}
\label{dfn:F and Ge}
{\rm 
We consider the $\mathbb{C}$-linear subspace $\mathcal{F}\subseteq \mathcal{O}((\mathbb{C}^*)^2)$ consisting of 
all $W$-invariant holomorphic functions $f(z)$ such that 
\begin{equation}
\label{eq:QP-F}
T_{p,z_1}f(z)=f(z)(pz_1^2z_2^3)^{-2}\quad\mbox{ and }\quad T_{p,z_2}f(z)=f(z)(p^3z_1^3z_2^6)^{-2}.
\end{equation}
We set $W_{0}=\la s_2, w_0\ra \subset W$, where $w_0=(s_1s_2)^3$ is the longest element of $W$. 
For $\epsilon \in\{-1,1\}$ we consider the $\mathbb{C}$-subspace $\mathcal{G}_{\epsilon}\subseteq z_{2}^{1\over 2}\mathcal{O}((\mathbb{C}^*)^2)$ consisting of 
all $W_0$-invariant functions $f(z)\in z_{2}^{1\over 2}\mathcal{O}((\mathbb{C}^*)^2)$ such that 
\begin{equation}
\label{eq:QP-G}
T_{p,z_1}f(z)=f(z)(pz_1^2z_2^3)^{-1}
\quad\mbox{ and }\quad 
T_{p,z_2}f(z)=f(z)(p^3z_1^3z_2^6)^{-1}(p^2z_2^4)^{-1}\epsilon.
\end{equation}
}
\end{dfn}
{\it Remark.} It can be verified directly that if $f(z)\in  \mathcal{M}((\mathbb{C}^*)^2)$ has the quasi-periodicity \eqref{eq:QP-F}, 
so does the function $w.f(z)$ for any $w\in W$.
Also if $f(z)\in z_{2}^{1\over 2}\mathcal{M}((\mathbb{C}^*)^2)$ satisfies \eqref{eq:QP-G}, 
so does $w.f(z)$ for any $w\in W_0$.

\begin{lem}
Under the condition $(a_1a_2a_3a_4a_5)q^{1\over 2}=\epsilon p^{1\over 2}$, 
if $\varphi(z)\in \mathcal{G}_{\epsilon}$, then $\nabla_{\rm sym}\varphi(z)\in \mathcal{F}$.
\end{lem}
{\it Proof.}  
Suppose that $\varphi(z)\in z_2^{1\over 2}\mathcal{M}((\mathbb{C}^*)^2)$ satisfies \eqref{eq:QP-G}. Then we have
\begin{equation*}
\begin{split}
T_{p,z_1}T_{q,z_1}^{1\over 2}\varphi(z)&=T_{q,z_1}^{1\over 2}\varphi(z)\, q^{-1}(pz_1^2z_2^3)^{-1},\\
T_{p,z_2}T_{q,z_1}^{1\over 2}\varphi(z)&=T_{q,z_1}^{1\over 2}\varphi(z)
\,q^{-{3\over 2}}(p^3z_1^3z_2^6)^{-1}(p^2z_2^4)^{-1}\epsilon.
\end{split}
\end{equation*}
Since $f^+(z)$ satisfies \eqref{eq:QP-f+}, we have 
\begin{equation*}
\begin{split}
T_{p,z_1}\big(f^+(z)T_{q,z_1}^{1\over 2}\varphi(z)\big)
&=f^+(z)T_{q,z_1}^{1\over 2}\varphi(z)\, q^{-1}p(a_1\cdots a_5 )^{-2}
(pz_1^2z_2^3)^{-2},\\
T_{p,z_2}\big(f^+(z)T_{q,z_1}^{1\over 2}\varphi(z)\big)
&=f^+(z)T_{q,z_1}^{1\over 2}\varphi(z)
\,\epsilon q^{-{3\over 2}}p^{3\over 2}(a_1\cdots a_5 )^{-3}
(p^3z_1^3z_2^6)^{-2}.
\end{split}
\end{equation*}
Similarly, from \eqref{eq:QP-f-} we have 
\begin{equation*}
\begin{split}
T_{p,z_1}\big(f^-(z)T_{q,z_1}^{-{1\over 2}}\varphi(z)\big)
&=f^-(z)T_{q,z_1}^{-{1\over 2}}\varphi(z)\, qp^{-1}(a_1\cdots a_5 )^{2}
(pz_1^2z_2^3)^{-2},\\
T_{p,z_2}\big(f^-(z)T_{q,z_1}^{-{1\over 2}}\varphi(z)\big)
&=f^-(z)T_{q,z_1}^{-{1\over 2}}\varphi(z)
\,\epsilon q^{3\over 2}p^{-{3\over 2}}(a_1\cdots a_5 )^{3}
(p^3z_1^3z_2^6)^{-2}.
\end{split}
\end{equation*}
Hence, under the condition $(a_1a_2a_3a_4a_5)q^{1\over 2}=\epsilon p^{1\over 2}$, 
both $f^+(z)T_{q,z_1}^{1\over 2}\varphi(z)$ and 
$f^-(z)T_{q,z_1}^{-{1\over 2}}\varphi(z)$ satisfy the same quasi-periodicity condition \eqref{eq:QP-F}, 
as well as $\nabla \varphi(z)=f^+(z)T_{q,z_1}^{1\over 2}\varphi(z)+f^-(z)T_{q,z_1}^{-{1\over 2}}\varphi(z)$. 
This implies that $\nabla_{\rm sym} \varphi(z)=\sum_{w\in W} w.(\nabla \varphi(z))$ satisfies 
\eqref{eq:QP-F}, since this quasi-periodicity is preserved by the action of $W$. 

We introduce the elliptic version of the Weyl denominator 
\begin{equation}
\label{eq:Delta(z;p)}
\Delta(z)=\Delta(z;p)=z_1^{-3}z_2^{-5}\theta(z_2,z_1z_2,z_1z_2^2,z_1,z_1z_2^3,z_1^2z_2^3;p), 
\end{equation}
which satisfies 
$
w.\Delta(z)={\rm sgn} (w)\Delta(z)
$ for all $w\in W$. Note that 
$$
\nabla_{\rm sym} \varphi(z)
=\frac{1}{\Delta(z)}
\sum_{w\in W}{\rm sgn} (w)\,
w.(\Delta(z)f^+(z)T_{q,z_1}^{1\over 2}\varphi(z)+\Delta(z)f^-(z)T_{q,z_1}^{-{1\over 2}}\varphi(z)).
$$
When $\varphi(z)$ belongs to $z_2^{1\over 2}\mathcal{O}((\mathbb{C}^*)^2)$, 
the numerator of the right-hand side is a quasi-periodic holomorphic function on $(\mathbb{C}^*)^2$. 
Since it is alternating with respect to the action of $W$, it is divisible by $\Delta(z)$. This means that 
$\nabla_{\rm sym} \varphi(z)$ is holomorphic on $(\mathbb{C}^*)^2$ and belongs to the space $\mathcal{F}$. \qed
\medskip
\begin{lem}
For $\varphi(z)\in \mathcal{G}_\epsilon$, $\nabla_{\rm sym}\varphi(z)$ is expressed as 
\begin{equation}
\label{eq:nabla_sym2}
\nabla_{\rm sym}\varphi(z)=4\sum_{k=0}^5f_k(z)\varphi_k(z),
\end{equation}
where $f_k(z)$ and $\varphi_k(z)$ 
are given by 
\begin{equation}
\label{eq:def fk(z)=}
f_k(z)=(s_1s_2)^k.f^{+}(z), \quad \varphi_k(z)=(s_1s_2)^k. T_{q,z_1}^{1\over 2}\varphi(z)
\end{equation}
for $k=0,1,\ldots,5$.
\end{lem}
{\it Proof.}
From the definition \eqref{eq:f-} of $f^-(z)$, we have 
\begin{equation}
\label{eq:f^-=f_3}
f^-(z)=f^+(z^{-1})=(s_1s_2)^3.f^+(z)=f_3(z). 
\end{equation}
Since $\varphi(z)\in \mathcal{G}_\epsilon$ is invariant under $w_0=(s_1s_2)^3$, i.e., $\varphi(z)=\varphi(z^{-1})$,
for $\varphi(z)\in \mathcal{G}_{\epsilon}$ we have 
\begin{equation}
\label{eq:T=varphi_3}
T_{q,z_1}^{-{1\over 2}}\varphi(z)=\varphi(q^{-{1\over 2}}z_1,z_2)
=\varphi(q^{{1\over 2}}z_1^{-1},z_2^{-1})
=\varphi_0(z^{-1})=(s_1s_2)^3.\varphi_0(z)=\varphi_3(z). 
\end{equation}
Applying \eqref{eq:f^-=f_3} and \eqref{eq:T=varphi_3} to the definition 
\eqref{eq:nabla} of $\nabla\varphi(z)$, for $\varphi(z)\in \mathcal{G}_{\epsilon}$ we have 
$$
\nabla\varphi(z)
=f^+(z)T_{q,z_1}^{{1\over 2}}\varphi(z)+f^+(z^{-1})T_{q,z_1}^{-{1\over 2}}\varphi(z)
=f_0(z)\varphi_0(z)+f_3(z)\varphi_3(z),
$$
so that we have 
$$
\nabla_{\rm sym}\varphi(z)
=\sum_{w\in W}w.(\nabla\varphi(z))=
\sum_{w\in W}w.(f_0(z)\varphi_0(z)+f_3(z)\varphi_3(z))
=2\sum_{w\in W}w.(f_0(z)\varphi_0(z)). 
$$
From the expression \eqref{eq:W} of $W$, this implies that 
\begin{equation}
\label{eq:nabla_sym3}
\nabla_{\rm sym}\varphi(z)=2\sum_{k=0}^5(s_1s_2)^k.
\Big(f_0(z)\varphi_0(z)+s_2.(f_0(z)\varphi_0(z))\Big).  
\end{equation}
Since $\varphi(z)\in \mathcal{G}_\epsilon$ is invariant under $s_2$, i.e., 
$\varphi(z_1,z_2)=\varphi(z_1z_2^3,z_2^{-1})$,
for $\varphi(z)\in \mathcal{G}_{\epsilon}$ we have 
\begin{equation}
\label{eq:s_2varphi_0}
s_2.\varphi_0(z)=\varphi(q^{{1\over 2}}z_1z_2^3,z_2^{-1})
=\varphi(q^{{1\over 2}}z_1,z_2)=\varphi_0(z).
\end{equation}
From \eqref{eq:s_2f+=f+}, \eqref{eq:s_2varphi_0} and \eqref{eq:nabla_sym3}, we therefore obtain 
$
\nabla_{\rm sym}\varphi(z)=4\sum_{k=0}^5(s_1s_2)^k.
(f_0(z)\varphi_0(z))  
$, 
which coincides with \eqref{eq:nabla_sym2}. \qed

\section{Proof of Lemma \ref{lem:qDE1}}\label{section:6}
The goal of this section is to give a proof of Lemma \ref{lem:qDE1} investigating 
the $\mathbb{C}$-linear mapping $\nabla_{\rm sym}: \mathcal{G}_\epsilon\to\mathcal{F}$ defined in the previous section. 
For that purpose we first clarify the structure of the target space $\mathcal{F}$ in Definition \ref{dfn:F and Ge}.  
\begin{lem}
\label{lem:dim F le 4}
$\dim_{\mathbb{C}}\mathcal{F}\le 4$.
\end{lem}
{\it Proof.} For arbitrary $f(z)\in\mathcal{F}$, since $f(z)$ is a holomorphic function of $z=(z_1,z_2)\in (\mathbb{C}^*)^2$, 
$f(z)$ can be expanded as Laurent series $\sum_{\lambda\in P}c_\lambda z^{\lambda}$, 
where $\lambda=\lambda_1\alpha_1+\lambda_2\alpha_2\in Q=P$. 
\begin{figure}[htbp]
 \begin{center}
\includegraphics[width=300pt]{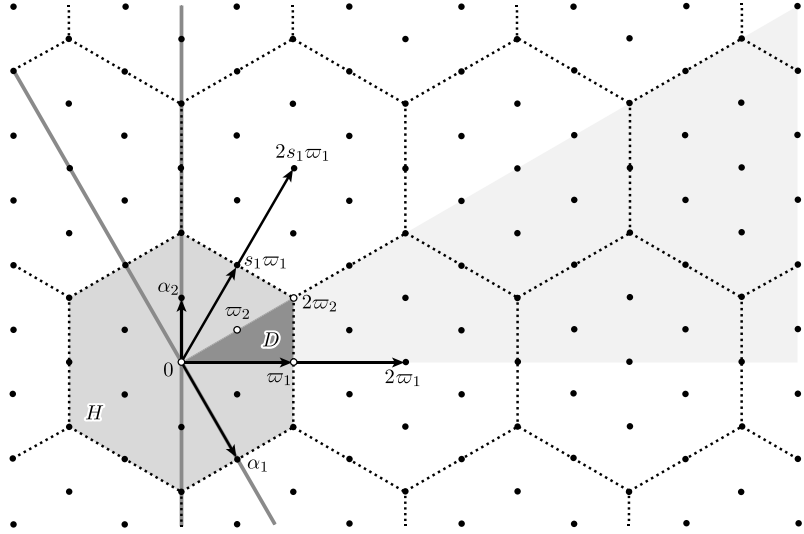}
\end{center}
\vspace{-10pt}
\caption{$D$, $H$ and the lattice points in $P$}
\label{fig:DHP}
\end{figure}
Let $D$ be the set defined by 
$D=\{\lambda\in P\,|\, (\lambda,\alpha_1)\ge 0,(\lambda,\alpha_2)\ge 0, (\lambda-\varpi_1,\varpi_1)\le 0\}$, 
which is the set of points in the triangle area consisting of three lines $(\lambda,\alpha_1)=0$, 
$(\lambda,\alpha_2)=0$ and $(\lambda-\varpi_1,\varpi_1)=0$.
Then we can immediately confirm that $D=\{0,\varpi_1,\varpi_2,2\varpi_2\}$. 
(See Figure \ref{fig:DHP}.)
Let $H$ be the $W$-orbit of $D$, i.e., $H=\bigcup_{w\in W}w D$, 
which is the set of points in the regular hexagon area consisting of six lines 
$(\lambda-(s_2s_1)^k\varpi_1,(s_2s_1)^k\varpi_1)=0$ for $k=0,1,\ldots,5$.
Since $f(z)\in \mathcal{F}$ is $W$-symmetric, i.e., $wf(z)=f(z)$ for $w\in W$,  
we have $c_{w\lambda}=c_\lambda$ for $w\in W$. 
This means that the coefficients $c_\lambda$ for $\lambda\in H$ are determined by $c_\lambda$ for $\lambda\in D$. 
On the other hand, since $f(z)\in \mathcal{F}$ satisfies that 
\begin{equation}
\label{eq:QP of F}
T_{p,z_1}f(z)=f(z)(pz^{\varpi_1})^{-2},
\end{equation}
we have 
$
\sum_{\lambda\in P}c_\lambda\,p^{\lambda_1}z^\lambda
=\sum_{\lambda\in P}c_\lambda\,p^{-2}z^{\lambda-2\varpi_1}. 
$
Equating the coefficients of $z_\lambda$ on both sides, we have $c_\lambda\,p^{\lambda_1}=c_{\lambda+2\varpi_1}p^{-2}$, i.e., 
\begin{equation}
\label{eq:coefficient of f--1}
c_{\lambda+2\varpi_1}=c_\lambda\,p^{\lambda_1+2} \quad (\lambda\in P). 
\end{equation}
Moreover, applying \eqref{eq:QP of F} to $f(z)=s_1.f(z)=\sum_{\lambda\in P}c_{s_1\lambda}z^\lambda$, we also have $c_{s_1\lambda}\,p^{\lambda_1}=c_{s_1(\lambda+2\varpi_1)}p^{-2}$, i.e., 
\begin{equation}
\label{eq:coefficient of f--2}
c_{\lambda+2(s_1\varpi_1)}=c_\lambda\,p^{\lambda_2-\lambda_1+2} \quad (\lambda\in P), 
\end{equation}
where $s_1\varpi_1=\alpha_1+3\alpha_2$. 
Combining \eqref{eq:coefficient of f--1} and \eqref{eq:coefficient of f--2}, 
for arbitrary $n_1,n_2\in \mathbb{Z}$ we have 
$$
c_{\lambda+2n_1\varpi_1+2n_2(s_1\varpi_1)}=c_\lambda\, p^{n_1(\lambda_1+2)+n_2(\lambda_2-\lambda_1+2)} \quad (\lambda\in P),
$$
so that the coefficients $c_\lambda$ for $\lambda\in H_{n_1,n_2}=\{\lambda+2n_1\varpi_1+2n_2(s_1\varpi_1)\,|\,\lambda\in H\}$ 
are determined by $c_\lambda$ for $\lambda\in H$. 
Since the lattice $P$ is covered by the sets $H_{n_1,n_2}$ for $(n_1,n_2)\in \mathbb{Z}^2$, i.e., 
$P=\bigcup_{(n_1,n_2)\in \mathbb{Z}^2}H_{n_1,n_2}$, 
the coefficients $c_\lambda$ for $\lambda\in P$ are also determined by $c_\lambda$ for $H$, and  
consequently, determined by $c_\lambda$ for $\lambda\in D$. 
Therefore we obtain $\dim_{\mathbb{C}}\mathcal{F}\le |D|=4$.
\qed
\par
\medskip
\noindent
{\it Remark.} In the next proposition, we prove that $\dim_{\mathbb{C}}\mathcal{F}=|D|=4$. 
In general, we can consider  
the $\mathbb{C}$-linear subspace $\mathcal{F}_n\subseteq \mathcal{O}((\mathbb{C}^*)^2)$ consisting of 
all $W$-invariant holomorphic functions $f(z)$ such that 
\begin{equation*}
\label{eq:QP-Fn}
T_{p,z_1}f(z)=f(z)(pz_1^2z_2^3)^{-n}\quad\mbox{ and }\quad T_{p,z_2}f(z)=f(z)(p^3z_1^3z_2^6)^{-n}. 
\end{equation*}
In the same way as above, it is actually confirmed that 
$$\dim_{\mathbb{C}}\mathcal{F}_n=|D_n|=
\begin{cases}
m(m+1)&(n=2m-1)\\
(m+1)^2&(n=2m)
\end{cases},$$
where 
$$D_n=\{\lambda\in P\,|\, (\lambda,\alpha_1)\ge 0,(\lambda,\alpha_2)\ge 0, (\lambda-\frac{n}{2}\varpi_1,\frac{n}{2}\varpi_1)\le 0\}.$$

\begin{prop}
\label{prop: basisF}
Let $F_k(z)$ be functions defined in \eqref{eq:F_k(z)}. 
For generic $a_1,a_2,a_3,a_4\in \mathbb{C}^*$, the set $\{F_1(z),\ldots,F_4(z)\}$ is a $\mathbb{C}$-basis of the space $\mathcal{F}$, i.e.,  
$$
\mathcal{F}=\mathbb{C}F_1(z)\oplus\mathbb{C}F_2(z)\oplus\mathbb{C}F_3(z)\oplus\mathbb{C}F_4(z). 
$$
In particular, $\dim_{\mathbb{C}}\mathcal{F}=4$.
\end{prop}

In the rest of this section, we omit the base $p$ in the notation $e(u,v;p)$ defined in \eqref{eq:e(u,v)}, so that 
$$
F_k(z)
=e(a_k,z_2)e(a_k,z_1z_2)e(a_k,z_1z_2^2). 
$$
In order to prove Proposition \ref{prop: basisF}, 
we define the points $\mathrm{p}_{ij}$ $(1\le i<j\le 5)$ in $(\mathbb{C}^*)^2$ given by
\begin{equation}
\label{eq:pij}
\mathrm{p}_{ij}=(a_i/a_j, a_j),
\end{equation}
i.e., $z_2=a_j$,  $z_1z_2=a_i$, $z_1z_2^2=a_ia_j$ if $z={\rm p}_{ij}$. 
By definition 
the values $F_k(z)$ at $z=\mathrm{p}_{ij}$ are given by
\begin{equation}
\begin{split}
\label{eq: Fk(p)}
F_k(\mathrm{p}_{ij})
&=e(a_k,a_i)e(a_k,a_j)e(a_k,a_ia_j)\\
&=a_k^{-3}\theta(a_ka_i,a_ka_i^{-1},a_ka_j,a_ka_j^{-1}, a_ka_ia_j,a_ka_i^{-1}a_j^{-1};p),
\end{split}
\end{equation}
and hence
\begin{equation}
\label{eq: Fk(p)=0}
k\in \{i,j\}\ \Longrightarrow \ F_k(\mathrm{p}_{ij})=0.
\end{equation}
This already implies that $F_k(z)$ $(k=1,2,3)$ are linearly independent for generic $a_1,a_2,a_3$. 
In fact, their values at $z=\mathrm{p}_{23},\mathrm{p}_{13},\mathrm{p}_{12}$ are given as follows.
\[
\begin{array}{c||ccc}
  & \mathrm{p}_{23}  & \mathrm{p}_{13}  & \mathrm{p}_{12}\\[3pt]
 \hline\\[-10pt]
F_1  &  * &  0 & 0\\[3pt]
F_2  &  0 &  * & 0\\[3pt]
F_3  &  0 & 0  & *\\[3pt]
\end{array}
\]
Let $G(z)$ be function defined by
\begin{equation}
\label{eq:def of G(z)}
G(z)=\frac{1}{\prod_{i=1}^3e(a_4,a_i)}
\left(F_4(z)-\frac{F_4(\mathrm{p}_{23})}{F_1(\mathrm{p}_{23})}F_1(z)
-\frac{F_4(\mathrm{p}_{13})}{F_2(\mathrm{p}_{13})}F_2(z)
-\frac{F_4(\mathrm{p}_{12})}{F_3(\mathrm{p}_{12})}F_3(z)
\right). 
\end{equation}
\begin{lem}
\label{lem:G(p)}
For $z=\mathrm{p}_{ij}\in (\mathbb{C}^*)^2$ $(1\le i<j\le 3)$ the function $G(z)$ satisfies 
\begin{equation}
\label{eq:G(p)=0}
G(\mathrm{p}_{12})=G(\mathrm{p}_{13})=G(\mathrm{p}_{23})=0.
\end{equation}
Moreover for $\{i,j,k\}=\{1,2,3\}$ 
and $x\in \mathbb{C}^*,$ if $z=(a_k/x,x),$ then we have
\begin{equation}
\label{eq:G(ak/x,x)}
G(a_k/x,x)=\frac{\theta(a_k^2,a_ka_ix,a_i/x,a_ka_jx,a_j/x;p)}
{\theta(a_k,a_ia_j,a_ia_ja_k,a_ia_k/a_j,a_ja_k/a_i;p)},
\end{equation}
which is independent of $a_4$. In particular, we have 
\begin{equation}
\label{eq:G(pk4)}
G(\mathrm{p}_{k4})=\frac{\theta(a_k^2,a_ka_ia_4,a_i/a_4,a_ka_ja_4,a_j/a_4;p)}{\theta(a_k,a_ia_j,a_ia_ja_k,a_ia_k/a_j,a_ja_k/a_i;p)}.
\end{equation}
\end{lem}
{\it Proof.} From \eqref{eq: Fk(p)=0} and \eqref{eq:def of G(z)}, we have \eqref{eq:G(p)=0}. 
Without loss of generality we prove \eqref{eq:G(ak/x,x)} for $k=1$. 
Using \eqref{eq: Fk(p)} we also have 
\begin{align}
\label{eq:G(a1/x,x)}
&G(a_1/x,x)
=
\frac{1}{\prod_{i=1}^3e(a_4,a_i)}
\left(F_4(a_1/x,x)
-\frac{F_4(\mathrm{p}_{13})}{F_2(\mathrm{p}_{13})}F_2(a_1/x,x)
-\frac{F_4(\mathrm{p}_{12})}{F_3(\mathrm{p}_{12})}F_3(a_1/x,x)
\right)\nonumber\\
&=
\frac{e(a_4,x)e(a_4,a_1x)}{e(a_4,a_2)e(a_4,a_3)}
-\frac{e(a_2,x)e(a_4,a_1a_3)e(a_2,a_1x)}{e(a_4,a_2)e(a_2,a_3)e(a_2,a_1a_3)}
-\frac{e(a_3,x)e(a_4,a_1a_2)e(a_3,a_1x)}{e(a_4,a_3)e(a_3,a_2)e(a_3,a_1a_2)}.
\end{align}
Denoting $g(x)=G(a_1/x,x)$ as a function of $x$,  
from \eqref{eq:G(a1/x,x)}
we immediately see that $g(x)$ is holomorphic on $\mathbb{C}^*$
and satisfies 
$$
g(px)=g(x)(p^2a_1^2x^4)^{-1}
$$
and $$g(a_1^{-1}a_2^{-1})=g(a_1^{-1}a_3^{-1})=g(a_2)=g(a_3)=0.$$ 
This implies that $g(x)$ is divisible by $\theta(a_1a_2x,a_1a_3x,a_2/x,a_3/x;p)$, 
i.e.,
\begin{equation}
\label{eq:g(x)=c theta}
g(x)=c\,\theta(a_1a_2x,a_1a_3x,a_2/x,a_3/x;p),
\end{equation}
where $c$ is a constant independent of $x$. We evaluate $g(a_2/a_1)$ in two ways. 
From \eqref{eq:G(a1/x,x)} we have  
\begin{align}
&g(a_2/a_1)=G(a_1^2/a_2,a_2/a_1)=
\frac{e(a_4,a_2/a_1)}{e(a_4,a_3)}
-\frac{e(a_3,a_2/a_1)e(a_4,a_1a_2)}{e(a_4,a_3)e(a_3,a_1a_2)}
\nonumber\\
&\quad=\frac{e(a_4,a_2/a_1)e(a_3,a_1a_2)-e(a_3,a_2/a_1)e(a_4,a_1a_2)}
{e(a_4,a_3)e(a_3,a_1a_2)}
=\frac{e(a_4,a_3)e(a_2/a_1,a_1a_2)}{e(a_4,a_3)e(a_3,a_1a_2)}
\quad\mbox{(using  \eqref{eq:R/W-relation})}
\nonumber\\
&\quad=\frac{e(a_2/a_1,a_1a_2)}{e(a_3,a_1a_2)}
=\frac{\theta(a_1^2,a_2^2;p)}{\theta(a_1a_2a_3,a_1a_2/a_3;p)}. 
\label{eq:G(p12*)00}
\end{align} 
On the other hand, \eqref{eq:g(x)=c theta} implies 
$
g(a_2/a_1)=
c\,\theta(a_2^2,a_2a_3,a_1,a_1a_3/a_2;p),
$
so that we have 
$$c=\frac{\theta(a_1^2,a_2^2;p)/\theta(a_1a_2a_3,a_1a_2/a_3;p)}
{\theta(a_2^2,a_2a_3,a_1,a_1a_3/a_2;p)}
=\frac{\theta(a_1^2;p)}
{\theta(a_1,a_2a_3,a_1a_2a_3,a_1a_2/a_3,a_1a_3/a_2;p)}
.$$
We therefore obtain 
$$
g(x)=
G(a_1/x,x)=
\frac{\theta(a_1^2,a_1a_2x,a_1a_3x,a_2/x,a_3/x;p)}
{\theta(a_1,a_2a_3,a_1a_2a_3,a_1a_2/a_3,a_1a_3/a_2;p)},
$$
which coincides with \eqref{eq:G(ak/x,x)} for $k=1$.
\qed 
\medskip
\noindent
{\bf Proof of  Proposition \ref{prop: basisF}.}
From \eqref{eq:G(p)=0} and \eqref{eq:G(pk4)} of Lemma \ref{lem:G(p)}, we see that $F_1(z),F_2(z),F_3(z),G(z)$
are linearly independent for generic $a_1,a_2,a_3, a_4$.
\[
\begin{array}{c||cccc}
  & \mathrm{p}_{23}  & \mathrm{p}_{13}  & \mathrm{p}_{12}& \mathrm{p}_{14}\\[3pt]
 \hline\\[-10pt]
F_1  &  * &  0 & 0& 0\\[3pt]
F_2  &  0 &  * & 0& *\\[3pt]
F_3  &  0 & 0  & *& *\\[3pt]
G  &  0 & 0  & 0& *\\[3pt]
\end{array}
\]
Hence $F_1(z),F_2(z),F_3(z),F_4(z)$
are also linearly independent for generic $a_1,a_2,a_3, a_4$. 
From Lemma \ref{lem:dim F le 4}, we therefore obtain that $\{F_1(z),\ldots,F_4(z)\}$ is a $\mathbb{C}$-basis of $\mathcal{F}$. 
This completes the proof of Proposition \ref{prop: basisF}.
\qed
\par
\medskip
By the definition \eqref{eq:def of G(z)} of $G(z)$ and  Proposition \ref{prop: basisF},  
the set $\{F_1(z),F_2(z),F_3(z),G(z)\}$ is also a $\mathbb{C}$-basis of the space $\mathcal{F}$, i.e.,
$\mathcal{F}
=\mathbb{C}F_1(z)\oplus\mathbb{C}F_2(z)\oplus\mathbb{C}F_3(z)\oplus\mathbb{C}G(z)
$. 
Here we slightly deform $F_3(z)$ in the set $\{F_1(z),F_2(z),F_3(z),G(z)\}$ 
as follows. 
Let ${\rm p}_{ij}^*$ $(1\le i<j\le 5)$ be points in $(\mathbb{C}^*)^2$ given by 
\begin{equation}
\label{eq:pij*}
{\rm p}_{ij}^*=(a_i^2/a_j,a_j/a_i), 
\end{equation}
so that $z_2=a_j/a_i$,  $z_1z_2=a_i$, $z_1z_2^2=a_j$ if $z={\rm p}_{ij}^*$. 
We remark that $F_1({\rm p}_{12}^*)=F_2({\rm p}_{12}^*)=0$ and 
$G({\rm p}_{12}^*)$ is evaluated as 
\begin{equation}
\label{eq:G(p12*)}
G({\rm p}_{12}^*)
=\frac{\theta(a_1^2,a_2^2;p)}{\theta(a_1a_2a_3,a_1a_2/a_3;p)}, 
\end{equation}
which is already confirmed in \eqref{eq:G(p12*)00}. 
We now define $F'_3(z)$ as 
$$
F'_3(z)=F_3(z)-\frac{F_3({\rm p}_{12}^*)}{G({\rm p}_{12}^*)}G(z),
$$
which satisfies  
\begin{equation}
\label{eq:F'3(p12),F'3(p12*)}
F'_3({\rm p}_{12})=F_3({\rm p}_{12})=e(a_3,a_2/a_1)e(a_3,a_1)e(a_3,a_2),\quad 
F'_3({\rm p}_{12}^*)=0.
\end{equation}
\[
\begin{array}{c||cccc}
  & \mathrm{p}_{23}  & \mathrm{p}_{13}  & \mathrm{p}_{12}& \mathrm{p}_{12}^*\\[3pt]
 \hline\\[-10pt]
F_1  &  * &  0 & 0& 0\\[3pt]
F_2  &  0 &  * & 0& 0\\[3pt]
F'_3  &  0 & 0  & *& 0\\[3pt]
G  &  0 & 0  & 0& *\\[3pt]
\end{array}
\]
Therefore the set $\{F_1(z),F_2(z),F'_3(z),G(z)\}$ is also 
a $\mathbb{C}$-basis of the space $\mathcal{F}$, i.e.,
\begin{equation}
\label{eq:F=F1+F2+F'3+G}
\mathcal{F}
=\mathbb{C}F_1(z)\oplus\mathbb{C}F_2(z)\oplus\mathbb{C}F'_3(z)\oplus\mathbb{C}G(z), 
\end{equation}
which will be used in the proof of Lemma \ref{lem:nabla_sym phi_12}.
\par
\medskip
Let $\phi_{ij}(z)$ and $\phi'_{ij}(z)$ be functions defined by
\begin{align}
\phi_{ij}(z)&=z_2^{-{1\over 2}}\theta(-z_2,-\epsilon p^{1\over 2}z_1z_2,-\epsilon p^{1\over 2}z_1z_2^2;p)e(a_i,z_2)e(a_j,z_2),
\label{eq:phi_ij}\\
\phi'_{ij}(z)&=z_1^{-1}z_2^{-{3\over 2}}\theta(-\epsilon p^{1\over 2}z_2,-z_1z_2,-z_1z_2^2;p)e(a_i,z_2)e(a_j,z_2).
\label{eq:phi_ij*}
\end{align}
From direct calculation we can immediately confirm that 
\begin{equation}
\label{eq:phi in Gepsilon}
\phi_{ij}(z), \phi'_{ij}(z)\in\mathcal{G}_{\epsilon}. 
\end{equation}
\noindent
\begin{prop}
For generic $a_1,a_2,a_3$, the set 
$\{\phi_{12}(z),\phi_{13}(z),\phi_{23}(z),\phi'_{12}(z),\phi'_{13}(z),\phi'_{23}(z)\}$ is a $\mathbb{C}$-basis of the space $\mathcal{G}_{\epsilon}$. 
In particular, $\dim_{\mathbb{C}}\mathcal{G}_{\epsilon}=6$.
\end{prop}
{\it Remark.} We omit the proof of this proposition since we need the fact 
\eqref{eq:phi in Gepsilon} only  for the succeeding discussions.  
\begin{lem} 
\label{lem:nabla_sym phi_12}
Under the condition $(a_1a_2a_3a_4a_5)q^{1\over 2}=\epsilon p^{1\over 2}$, 
$\nabla_{\rm sym}\phi_{12}(z)$ and $\nabla_{\rm sym}\phi'_{12}(z)$ are expanded as 
\begin{align}
\frac{\nabla_{\rm sym}\phi_{12}(z)}{4}&=c_1F_1(z)+c_2F_2(z)+c_{12}G(z),
\label{nabla_sym phi_{12}(z)}
\\
\frac{\nabla_{\rm sym}\phi'_{12}(z)}{4}&=c_1'F_1(z)+c_2'F_2(z)+c_{12}'G(z),
\label{nabla_sym phi'_{12}(z)}
\end{align}
where the coefficients are given by 
\begin{align}
c_1&=
\theta(-a_3,-\epsilon  p^{1\over 2}q^{1\over 2}a_2,-\epsilon  p^{1\over 2}q^{1\over 2}a_2a_3;p)
\frac{\theta(a_2^2,a_2a_3a_4,a_2a_3a_5;p)}{\theta(a_2,a_2/a_1,a_2a_3/a_1;p)}
\prod_{k=3}^5\theta(a_2a_k;p),
\label{eq:c1}\\
c_2&=
\theta(-a_3,-\epsilon  p^{1\over 2}q^{1\over 2}a_1,-\epsilon  p^{1\over 2}q^{1\over 2}a_1a_3;p)
\frac{\theta(a_1^2,a_1a_3a_4,a_1a_3a_5;p)}{\theta(a_1,a_1/a_2,a_1a_3/a_2;p)}
\prod_{k=3}^5\theta(a_1a_k;p),
\label{eq:c2}\\
c_{12}&=a_2\theta(-a_1/a_2,-\epsilon p^{1\over 2}q^{1\over 2}a_1,-\epsilon p^{1\over 2}q^{1\over 2}a_2;p)
\frac{\theta(a_1a_2a_3,a_1a_2/a_3,a_1a_2;p)}{a_1^2a_2^2}
\prod_{k=3}^5\theta(a_1a_k,a_2a_k;p),
\label{eq:c12}
\\
c_1'&=
q^{-{1\over 2}}a_2^{-1}\theta(-\epsilon  p^{1\over 2}a_3,-q^{1\over 2}a_2,-q^{1\over 2}a_2a_3;p)
\frac{\theta(a_2^2,a_2a_3a_4,a_2a_3a_5;p)}{\theta(a_2,a_2/a_1,a_2a_3/a_1;p)}
\prod_{k=3}^5\theta(a_2a_k;p),
\label{eq:c1*}\\
c_2'&=
q^{-{1\over 2}}a_1^{-1}\theta(-\epsilon  p^{1\over 2}a_3,-q^{1\over 2}a_1,-q^{1\over 2}a_1a_3;p)
\frac{\theta(a_1^2,a_1a_3a_4,a_1a_3a_5;p)}{\theta(a_1,a_1/a_2,a_1a_3/a_2;p)}
\prod_{k=3}^5\theta(a_1a_k;p),
\label{eq:c2*}\\
c_{12}'&=q^{-{1\over 2}}
\theta(-\epsilon p^{1\over 2}a_1/a_2,-q^{1\over 2}a_1,-q^{1\over 2}a_2;p)
\frac{\theta(a_1a_2a_3,a_1a_2/a_3,a_1a_2;p)}{a_1^2a_2^2}
\prod_{k=3}^5\theta(a_1a_k,a_2a_k;p).
\label{eq:c12*}
\end{align}
\end{lem}
{\it Proof.}
For $\varphi(z)\in \mathcal{G}_{\epsilon}$ we simply use notation  
\begin{equation}
\label{eq:widetilde_varphi}
\widetilde\varphi(z)=\frac{\nabla_{\rm sym}\varphi(z)}{4}=\sum_{k=0}^5f_k(z)\varphi_k(z)\in \mathcal{F},
\end{equation}
where $f_k(z)$ and $\varphi_k(z)$ are given in \eqref{eq:def fk(z)=}. 
Our aim is to confirm that $
\widetilde\phi_{12}(z)$ and $\widetilde\phi'_{12}(z) 
$
are expanded as the right-hand side of 
\eqref{nabla_sym phi_{12}(z)} and \eqref{nabla_sym phi'_{12}(z)}, respectively.  
Before that 
we evaluate the special values 
$\widetilde\phi_{12}(\mathrm{p}_{1j})$, $\widetilde\phi_{12}'(\mathrm{p}_{1j})$ 
and 
$\widetilde\phi_{12}(\mathrm{p}_{12}^*)$, $\widetilde\phi'_{12}(\mathrm{p}_{12}^*)$.
\par
First we evaluate 
$\widetilde\phi_{12}(\mathrm{p}_{1j})$ and $\widetilde\phi_{12}'(\mathrm{p}_{1j})$. 
From \eqref{eq:pij}, if $z=\mathrm{p}_{ij}$, then $z_2=a_j$, $z_1z_2=a_i$. 
Hence, since $f_2(z)$, $f_3(z)$ have factors 
$
\theta(a_k z_1^{-1}z_2^{-1};p)
$ $(1\le k\le 5)$
and 
$f_4(z)$, $f_5(z)$ have factors 
$
\theta(a_k z_2^{-1};p)
$ $(1\le k\le 5)$, we have 
\begin{equation*}
f_2(\mathrm{p}_{ij})=f_3(\mathrm{p}_{ij})=f_4(\mathrm{p}_{ij})=f_5(\mathrm{p}_{ij})=0
\quad\mbox{for}\quad 1\le i<j\le 5. 
\end{equation*}
\par
\begin{equation*}
\begin{picture}(150,100)
(-25,0)
\put(50,50){\vector(2,1){30}}
\put(50,50){\vector(0,1){30}}
\put(50,50){\vector(-2,1){30}}
\put(50,50){\vector(-2,-1){30}}
\put(50,50){\vector(0,-1){30}}
\put(50,50){\vector(2,-1){30}}
\put(47,85){$z_2$}
\put(84,63){$z_1z_2^2$}
\put(84,32){$z_1z_2$}
\put(47,10){$z_2^{-1}$}
\put(-13,30){$z_1^{-1}\!z_2^{-2}$}
\put(-13,63){$z_1^{-1}\!z_2^{-1}$}
\put(87,2){$z_1$}
\put(84,94){$z_1z_2^3$}
\put(115,47){$z_1^2z_2^3$}
\put(0,2){$z_1^{-1}\!z_2^{-3}$}
\put(10,94){$z_1^{-1}$}
\put(-47,47){$z_1^{-2}z_2^{-3}$}
\put(72,48){$f_0$}
\put(58,68){$f_1$}
\put(32,68){$f_2$}
\put(22,48){$f_3$}
\put(32,28){$f_4$}
\put(58,28){$f_5$}
\end{picture}
\end{equation*}
From \eqref{eq:widetilde_varphi}, this implies that  
\begin{equation*}
\widetilde\varphi(\mathrm{p}_{ij})=f_0(\mathrm{p}_{ij})\varphi_0(\mathrm{p}_{ij})+
f_1(\mathrm{p}_{ij})\varphi_1(\mathrm{p}_{ij}).
\end{equation*}
By definition we have $\varphi_0(z)=\varphi(q^{1\over2}z_1,z_2)$ and  
$\varphi_1(z)=\varphi(q^{1\over2}z_1^2z_2^3,z_1^{-1}z_2^{-1})$.  
If we put $\varphi(z)=\phi_{12}(z)$, then the explicit forms of $\varphi_0(z)$ and $\varphi_1(z)$ are given as 
\begin{align*}
\varphi_0(z)
&=z_2^{-{1\over 2}}
\theta(-z_2,-\epsilon p^{1\over 2}q^{1\over 2}z_1z_2, -\epsilon p^{1\over 2}q^{1\over 2}z_1z_2^2;p)e(a_1,z_2)e(a_2,z_2),
\\
\varphi_1(z)
&=(z_1z_2)^{-{1\over 2}}\theta(-z_1z_2,-\epsilon p^{1\over 2}q^{1\over 2}z_2, -\epsilon p^{1\over 2}q^{1\over 2}z_1z_2^2;p)e(a_1,z_1z_2)e(a_2,z_1z_2),
\end{align*}
so that we have 
\begin{equation}
\begin{split}
\varphi_0(\mathrm{p}_{ij})
&=a_j^{-{1\over 2}}
\theta(-a_j,-\epsilon p^{1\over 2}q^{1\over 2}a_i, -\epsilon p^{1\over 2}q^{1\over 2}a_ia_j;p)e(a_1,a_j)e(a_2,a_j),
\\
\varphi_1(\mathrm{p}_{ij})
&=a_i^{-{1\over 2}}
\theta(-a_i,-\epsilon p^{1\over 2}q^{1\over 2}a_j, -\epsilon p^{1\over 2}q^{1\over 2}a_ia_j;p)e(a_1,a_i)e(a_2,a_i). 
\end{split}
\end{equation}
On the other hand, in the same way as above, if we put $\varphi(z)=\phi_{12}'(z)$, then we have
\begin{equation}
\begin{split}
\varphi_0(\mathrm{p}_{ij})
&=
q^{-{1\over 2}}
a_i^{-1}a_j^{-{1\over 2}}
\theta(-\epsilon p^{1\over 2}a_j,-q^{1\over 2}a_i, -q^{1\over 2}a_ia_j;p)e(a_1,a_j)e(a_2,a_j),
\\
\varphi_1(\mathrm{p}_{ij})
&=q^{-{1\over 2}}
a_j^{-1}a_i^{-{1\over 2}}
\theta(-\epsilon p^{1\over 2}a_i,-q^{1\over 2}a_j, -q^{1\over 2}a_ia_j;p)e(a_1,a_i)e(a_2,a_i). 
\end{split}
\end{equation}
Therefore,  
if $\varphi(z)=\phi_{12}(z),\phi_{12}'(z)$, then 
we eventually obtain 
$\widetilde\varphi(\mathrm{p}_{1j})=f_0(\mathrm{p}_{1j})\varphi_0(\mathrm{p}_{1j})$, i.e., 
\begin{equation}
\label{eq:widetilde_phi_{12}}
\begin{split}
\widetilde\phi_{12}(\mathrm{p}_{1j})
&=f^+(\mathrm{p}_{1j})\phi_{12}(q^{1\over 2}a_1/a_j,a_j),\\
\widetilde\phi_{12}'(\mathrm{p}_{1j})
&=f^+(\mathrm{p}_{1j})\phi_{12}'(q^{1\over 2}a_1/a_j,a_j)
\end{split}
\end{equation}
for $j=2,3,4,5$, in particular, $\widetilde\phi_{12}(\mathrm{p}_{12})=\widetilde\phi_{12}'(\mathrm{p}_{12})=0$.
\par
Next we evaluate 
$\widetilde\phi_{12}(\mathrm{p}_{12}^*)$ and $\widetilde\phi'_{12}(\mathrm{p}_{12}^*)$.
From \eqref{eq:pij*}, if $z=\mathrm{p}_{12}^*$, then $z_1z_2=a_1$, $z_1z_2^2=a_2$. 
Since $f_2(z)$, $f_3(z)$ have factors 
$
\theta(a_k z_1^{-1}z_2^{-1};p)
$ $(1\le k\le 5)$
and 
$f_3(z)$, $f_4(z)$ have factors 
$
\theta(a_k a_1^{-1}z_2^{-2};p)
$ $(1\le k\le 5)$, we have 
\begin{equation*}
f_2(\mathrm{p}_{12}^*)=f_3(\mathrm{p}_{12}^*)=f_4(\mathrm{p}_{12}^*)=0
\quad\mbox{for}\quad 1\le i<j\le 5. 
\end{equation*}
From \eqref{eq:widetilde_varphi}, this implies that  
\begin{equation*}
\widetilde\varphi(\mathrm{p}_{12}^*)
=f_0(\mathrm{p}_{12}^*)\varphi_0(\mathrm{p}_{12}^*)+
f_1(\mathrm{p}_{12}^*)\varphi_1(\mathrm{p}_{12}^*)+
f_5(\mathrm{p}_{12}^*)\varphi_5(\mathrm{p}_{12}^*).
\end{equation*}
When $\varphi(z)=\phi_{12}(z)$ or $\phi_{12}'(z)$ 
the functions $\varphi_0(z)$, $\varphi_1(z)$ and $\varphi_5(z)$
have factors $e(a_1,z_2)e(a_2,z_2)$, $e(a_1,z_1z_2)e(a_2,z_2z_2)$ and 
$e(a_1,z_1z_2^2)e(a_2,z_1z_2^2)$, respectively. This implies 
$\varphi_1(\mathrm{p}_{12}^*)=0$ and $\varphi_5(\mathrm{p}_{12}^*)=0$. 
We therefore obtain $\widetilde\varphi(\mathrm{p}_{12}^*)
=f_0(\mathrm{p}_{12}^*)\varphi_0(\mathrm{p}_{12}^*)$, i.e., 
\begin{equation}
\label{eq:widetilde_phi_{12}*}
\begin{split}
\widetilde\phi_{12}(\mathrm{p}_{12}^*)
&=f^+(\mathrm{p}_{12}^*)\phi_{12}(q^{1\over 2}a_1^2/a_2,a_2/a_1),\\
\widetilde\phi_{12}'(\mathrm{p}_{12}^*)
&=f^+(\mathrm{p}_{12}^*)\phi_{12}'(q^{1\over 2}a_1^2/a_2,a_2/a_1).
\end{split}
\end{equation}
\par
We now calculate the expansions of 
$
\widetilde\phi_{12}(z)$ and $\widetilde\phi'_{12}(z). 
$
From \eqref{eq:F=F1+F2+F'3+G}, $\widetilde\phi_{12}(z), \widetilde\phi'_{12}(z)\in \mathcal{F}$
are expanded in terms of $F_1(z),F_2(z),F'_3(z),G(z)$, i.e., 
\begin{align}
\widetilde\phi_{12}(z)&=c_1F_1(z)+c_2F_2(z)+c_3F'_3(z)+c_{12}G(z),
\label{eq:widetilde_phi_{12}_2}
\\
\widetilde\phi_{12}'(z)&=c_1'F_1(z)+c_2'F_2(z)+c_3'F'_3(z)+c'_{12}G(z),
\end{align}
where $c_1,c_2,c_3,c_{12},c'_1,c'_2,c'_3,c'_{12}$ are some constants independent of $z$. Our aim is to show that 
$c_3=c_3'=0$ and the others are given by 
\eqref{eq:c1}--\eqref{eq:c12*}.  
Applying \eqref{eq: Fk(p)=0}, \eqref{eq:G(p)=0} and \eqref{eq:F'3(p12),F'3(p12*)} 
to \eqref{eq:widetilde_phi_{12}_2}, we have the expressions 
$$
c_1=\frac{\widetilde\phi_{12}(\mathrm{p}_{23})}{F_1(\mathrm{p}_{23})},\quad 
c_2=\frac{\widetilde\phi_{12}(\mathrm{p}_{13})}{F_2(\mathrm{p}_{13})},\quad 
c_3=\frac{\widetilde\phi_{12}(\mathrm{p}_{12})}{F'_3(\mathrm{p}_{12})},\quad 
c_{12}=\frac{
\widetilde\phi_{12}(\mathrm{p}_{12}^*)}{G(\mathrm{p}_{12}^*)}.%
$$
In particular, we immediately have $c_3=0$ because $\widetilde\phi_{12}(\mathrm{p}_{12})=0$. 
Furthermore, from 
\eqref{eq: Fk(p)}
and 
\eqref{eq:widetilde_phi_{12}} we obtain 
\begin{align*}
c_2&=\frac{\widetilde\phi_{12}(\mathrm{p}_{13})}{F_2(\mathrm{p}_{13})}=
\frac{f^+(\mathrm{p}_{13})\phi_{12}(q^{1\over 2}a_1/a_3,a_3)}{F_2(\mathrm{p}_{13})}
\\
&=\frac{(a_1^2a_3)^{-{1\over 2}}\theta(a_1a_2a_3,a_1a_3a_4,a_1a_3a_5;p)
\prod_{k=1}^5\theta(a_1a_k;p)}{
e(a_2,a_1)e(a_2,a_3)e(a_2,a_1a_3)\theta(a_1,a_1a_3,a_1/a_3;p)},
\end{align*}
which is equal to \eqref{eq:c2}. 
Since $\widetilde\phi_{12}(z)$ is symmetric with respect to the transposition of $a_1$ and $a_2$, 
we obtain the expression \eqref{eq:c1} of $c_1$ from \eqref{eq:c2} of $c_2$ interchanging $a_1$ and $a_2$. 
From \eqref{eq:G(p12*)} and \eqref{eq:widetilde_phi_{12}*} 
we have 
\begin{align*}
c_{12}
&=\frac{
\widetilde\phi_{12}(\mathrm{p}_{12}^*)}{G(\mathrm{p}_{12}^*)}
=\frac{f^+(\mathrm{p}_{12}^*)\phi_{12}(q^{1\over 2}a_1^2/a_2,a_2/a_1)}{G(\mathrm{p}_{12}^*)}\\
&=\frac{(a_1a_2)^{-\frac{1}{2}}\prod_{k=1}^5\theta(a_1a_k,a_2a_k;p)}
{\theta(a_1,a_2,a_1^2/a_2,a_2^2/a_1,a_1a_2;p)}
\frac{\theta(a_1a_2a_3,a_1a_2/a_3;p)}{\theta(a_1^2,a_2^2;p)}\\
&\quad\times(a_2/a_1)^{-{1\over 2}}
\theta(-a_2/a_1,-\epsilon p^{1\over 2}q^{1\over 2}a_1, -\epsilon p^{1\over 2}q^{1\over 2}a_2;p)e(a_1,a_2/a_1)e(a_2,a_2/a_1),
\end{align*}
which is equal to \eqref{eq:c12}.
\par
On the other hand, for $c'_1$, $c'_2$ and $c'_{12}$ we shall omit the process to calculate the expressions \eqref{eq:c1*}, \eqref{eq:c2*} and \eqref{eq:c12*} 
since it is almost the same way as above for $c_1$, $c_2$ and $c_{12}$. 
\qed
\begin{prop}
\label{prop:F1=C1G}
Under the condition $(a_1a_2a_3a_4a_5)^2q=p$, for $k=1,2,3$ we have 
\begin{equation}
\label{eq:F1 eqiv C1G}
F_k(z)\equiv C_k G(z) \quad ({\rm mod}\ \nabla_{\rm sym}\mathcal{G}_\epsilon), 
\end{equation}
where the constant $C_k$ is given as
\begin{equation*}
C_k=\frac{\theta(qa_k^2,a_i,a_j,a_ka_ia_j,a_ka_4a_5,a_ia_j/a_k,a_ka_j/a_i,a_ka_i/a_j;p)}{a_k^3\theta(qa_1a_2a_3,a_i^2,a_j^2,a_ia_ja_4,a_ia_ja_5,a_ia_ja_4a_5;p)}
\prod_{\substack{1\le l\le 5\\[1pt] l\ne k}}
\theta(a_ka_l;p),
\end{equation*}
when $\{i,j,k\}=\{1,2,3\}$. In particular, we have 
\begin{equation}
\label{eq:F1=C1G}
\la F_k(z) \ra
=C_k\la G(z) \ra \quad (k=1,2,3). 
\end{equation}
\end{prop}
\noindent{\it Proof.} 
Since \eqref{eq:F1=C1G} follows from \eqref{eq:F1 eqiv C1G} by  Lemma \ref{lem:nabla=0}, 
we show \eqref{eq:F1 eqiv C1G}. 
By the symmetry with respect to the indices $1,2,3$, it suffices to consider the case where $k=1$ with 
\begin{equation}
\label{eq:C1}
C_1=\frac{\theta(qa_1^2,a_1a_2,a_2,a_3,a_1a_2a_3,a_1a_4a_5,a_2a_3/a_1,a_1a_3/a_2,a_1a_2/a_3;p)}{a_1^3\theta(qa_1a_2a_3,a_2^2,a_3^2,a_2a_3a_4,a_2a_3a_5,a_2a_3a_4a_5;p)}
\prod_{i=3}^5\theta(a_1a_i;p).
\end{equation}
From the condition $(a_1a_2a_3a_4a_5)^2q=p$, 
we suppose that $(a_1a_2a_3a_4a_5)q^{1\over 2}=\epsilon p^{1\over 2}$, where $\epsilon \in\{-1,1\}$. 
Then, by  Lemma \ref{lem:nabla_sym phi_12}, 
we have 
$$
c_1 F_1(z)+c_2F_2(z)+c_{12}G(z)\equiv 0\quad ({\rm mod}\ \nabla_{\rm sym}\mathcal{G}_\epsilon),
$$
$$
c'_1F_1(z)+c'_2F_2(z)+c'_{12} G(z)\equiv 0\quad ({\rm mod}\ \nabla_{\rm sym}\mathcal{G}_\epsilon),
$$
where $c_1, c_2, c_{12}$, $c'_1, c'_2, c'_{12}$ are explicitly given by 
\eqref{eq:c1}--\eqref{eq:c12*}. 
Eliminating $F_2(z)$ from the above equations, we have 
$$F_1(z) \equiv \frac{c'_{12}c_2-c_{12}c'_2}{c_1c'_2-c'_1c_2} G(z)
\quad ({\rm mod}\ \nabla_{\rm sym}\mathcal{G}_\epsilon),$$
since one can verify $c_1c'_2-c'_1c_2\ne 0$ in the course of computation.  
Therefore $C_1$ in \eqref{eq:F1=C1G} is calculated as 
\begin{align}
\label{eq:C1=...N/D}
C_1
&=\frac{c'_{12}c_2-c_{12}c'_2}{c_1c'_2-c'_1c_2}
\nonumber\\
&=
\frac{\theta(a_1/a_2,a_1a_2a_3,a_1a_2/a_3,a_2,a_2a_3/a_1,-q^{1\over 2}a_1,-\epsilon p^{1\over 2}q^{1\over 2}a_1;p)\prod_{i=2}^5\theta(a_1a_i;p)}
{a_1^3\,\theta(a_2^2,a_2a_3a_4,a_2a_3a_5,-a_3,-\epsilon p^{1\over 2}a_3;p)}
\frac{N}{D},
\end{align}
where $N$ and $D$ are expressed as 
\begin{align}
N&=
a_2
\theta(-a_1/a_2,
-\epsilon p^{1\over 2}q^{1\over 2}a_2,-\epsilon p^{1\over 2}a_3,-q^{1\over 2}a_1a_3;p)
\nonumber\\
&\hspace{30pt}
-a_1
\theta(-\epsilon p^{1\over 2}a_1/a_2,
-q^{1\over 2}a_2,-a_3,-\epsilon p^{1\over 2}q^{1\over 2}a_1a_3;p),
\label{eq:N-1}\\
D&=
a_2\theta(-\epsilon p^{1\over 2}q^{1\over 2}a_2,-\epsilon p^{1\over 2}q^{1\over 2}a_2a_3,-q^{1\over 2}a_1,-q^{1\over 2}a_1a_3;p)
\nonumber\\
&\hspace{30pt}-
a_1\theta(-q^{1\over 2}a_2,-q^{1\over 2}a_2a_3,-\epsilon p^{1\over 2}q^{1\over 2}a_1,-\epsilon p^{1\over 2}q^{1\over 2}a_1a_3;p).
\label{eq:D-1}
\end{align}
Here we can confirm that $N$ and $D$ are factorized as  
\begin{align}
N&=a_2\theta(\epsilon p^{1\over 2},\epsilon p^{1\over 2}q^{1\over 2}a_2a_3,q^{1\over 2}a_1, a_1a_3/a_2;p),
\label{eq:N-2}
\\
D&=a_2\theta(\epsilon p^{1\over 2}, \epsilon p^{1\over 2}a_3, a_1/a_2, qa_1a_2a_3;p),
\label{eq:D-2}
\end{align}
respectively. 
In fact, \eqref{eq:N-2} is confirmed as follows. 
Denoting by $g(u)$ the right-hand side of \eqref{eq:N-1} as a function of $u=a_1$,  
$g(u)$ satisfies 
$$
g(pu)=g(u)(q^{1\over 2}a_3/a_2)^{-1}(u^2)^{-1}
$$
and $g(q^{-{1\over 2}})=0$. This implies that 
$g(u)$ is divisible by $\theta(q^{{1\over 2}}u,a_3u/a_2;p)$, i.e.,
$$
g(u)=c\,\theta(q^{{1\over 2}}u,a_3u/a_2;p),
$$
where $c$ is a constant independent of $u$. 
Evaluating $g(-\epsilon p^{1\over 2}a_2)$ in two ways, we have  
$$
a_2
\theta(\epsilon p^{1\over 2},
-\epsilon p^{1\over 2}q^{1\over 2}a_2,-\epsilon p^{1\over 2}a_3,\epsilon p^{1\over 2}q^{1\over 2}a_2a_3;p)
=c\,\theta(-\epsilon p^{1\over 2}q^{1\over 2}a_2,-\epsilon p^{1\over 2}a_3;p),
$$
so that we have 
$c=a_2
\theta(\epsilon p^{1\over 2},\epsilon p^{1\over 2}q^{1\over 2}a_2a_3;p).
$
We therefore obtain 
$$
g(u)=a_2
\theta(\epsilon p^{1\over 2},\epsilon p^{1\over 2}q^{1\over 2}a_2a_3,q^{{1\over 2}}u,a_3u/a_2;p),
$$
which is equivalent to the right-hand side of \eqref{eq:N-2}.
For the case \eqref{eq:D-2} of $D$, we can actually deduce it from \eqref{eq:D-1} 
in the same way as \eqref{eq:N-2} of $N$ above.
\par 
From \eqref{eq:C1=...N/D}, \eqref{eq:N-2} and \eqref{eq:D-2}, we consequently obtain 
\begin{align*}
C_1
&=
\frac{\theta(a_1a_2a_3,a_1a_2/a_3,a_2,a_1/a_2,a_2a_3/a_1,-q^{1\over 2}a_1,-\epsilon p^{1\over 2}q^{1\over 2}a_1;p)\prod_{i=2}^5\theta(a_1a_i;p)}
{a_1^3\,\theta(a_2^2,a_2a_3a_4,a_2a_3a_5,-a_3,-\epsilon p^{1\over 2}a_3;p)}
\\
&\hspace{20pt}\times \frac{\theta(\epsilon p^{1\over 2}q^{1\over 2}a_2a_3,q^{1\over 2}a_1, a_1a_3/a_2;p)}
{\theta( \epsilon p^{1\over 2}a_3, a_1/a_2, qa_1a_2a_3;p)}\\
&=\frac{\theta(a_1a_2a_3,a_1a_2/a_3,a_2,a_2a_3/a_1,-q^{1\over 2}a_1,
-\epsilon p^{1\over 2}q^{1\over 2}a_1,
\epsilon p^{1\over 2}q^{1\over 2}a_2a_3,q^{1\over 2}a_1, a_1a_3/a_2;p)}
{a_1^3\,\theta(a_2^2,a_2a_3a_4,a_2a_3a_5,-a_3,-\epsilon p^{1\over 2}a_3,
\epsilon p^{1\over 2}a_3, qa_1a_2a_3;p)}
\prod_{i=2}^5\theta(a_1a_i;p),
\end{align*}
which coincides with \eqref{eq:C1} by using the identities 
\begin{align*}
\theta(qa_1^2;p)&=\theta(q^{1\over 2}a_1,-q^{1\over 2}a_1,\epsilon p^{1\over 2}q^{1\over 2}a_1,-\epsilon p^{1\over 2}q^{1\over 2}a_1;p),\\
\theta(a_3^2;p)&=\theta(a_3,-a_3,\epsilon p^{1\over 2}a_3,-\epsilon p^{1\over 2}a_3;p)
\end{align*}
and the relations 
\begin{align*}
\theta(\epsilon p^{1\over 2}q^{1\over 2}a_1;p)
&=\theta(\epsilon p^{1\over 2}q^{-{1\over 2}}a_1^{-1};p)=
\theta(a_2a_3a_4a_5;p),\\
\theta(\epsilon p^{1\over 2}q^{1\over 2}a_2a_3;p)
&=\theta(\epsilon p^{1\over 2}q^{-{1\over 2}}a_2^{-1}a_3^{-1};p)
=\theta(a_1a_4a_5;p)
\end{align*}
under the condition $(a_1a_2a_3a_4a_5)q^{1\over 2}=\epsilon p^{1\over 2}$. 
Thus we finally obtain the expression \eqref{eq:C1} of $C_1$ independent of $\epsilon$. 
\qed
\par
\medskip
\noindent{\bf Proof of Lemma \ref{lem:qDE1}.} 
By Proposition \ref{prop:F1=C1G}, 
under the condition $(a_1a_2a_3a_4a_5)^2q=p$, 
we have the relation between $\la F_2(z) \ra$ and $\la G(z)\ra$ as
$
\la F_2(z) \ra
=C_2\la G(z) \ra
$, 
where 
\begin{equation}
\label{eq:C2}
C_2=\frac{\theta(qa_2^2,a_1a_2,a_1,a_3,a_1a_2a_3,a_2a_4a_5,a_1a_3/a_2,a_2a_3/a_1,a_1a_2/a_3;p)}{a_2^3\,\theta(qa_1a_2a_3,a_1^2,a_3^2,a_1a_3a_4,a_2a_3a_5,a_1a_3a_4a_5;p)}
\prod_{i=3}^5\theta(a_2a_i;p).
\end{equation}
Therefore we obtain 
$\la F_1(z) \ra=(C_1/C_2)\la F_2(z)\ra$ and from \eqref{eq:C1} and \eqref{eq:C2}
 the coefficient $C_1/C_2$ coincides with that of \eqref{eq:F1/F2}
in Lemma \ref{lem:qDE1}. \qed
\par
\medskip 
\noindent
{\it Remark. } 
The quotient space 
$\mathcal{H}=\mathcal{F}/\nabla_{\rm sym}\mathcal{G}_\epsilon$  
can be interpreted as the 
$q$-difference de Rham cohomology associated with our integral $I(a)$. 
Note that Proposition \ref{prop:F1=C1G} means $\dim_{\mathbb{C}} \mathcal{H}\le 1$. 
Although one can directly verify that $\dim_{\mathbb{C}} \mathcal{H}=1$, we omit the proof since 
Theorem \ref{thm:eG2} should eventually imply $\mathcal{H}\ne 0$.

\section{\bf\boldmath Computation of the constant $b$}\label{section:7}
In this section we use the double-sign symbol like $\Gamma(uv^{\pm 1};p,q)=\Gamma(uv, uv^{-1};p,q)$
for abbreviation. 
\par
As before, we assume that the parameters satisfy the 
balancing condition $(a_1\cdots a_5)^2=pq$, 
and regard $a_{5}=\epsilon p^{1\over 2}q^{1\over 2}/a_1\cdots a_4$,
where $\epsilon \in \{-1,1\}$, 
as a function of $(a_1,\ldots,a_4)$.
By Theorem \ref{thm:I=cJ} 
we showed that 
the meromorphic functions 
$I(a_1,\ldots,a_5)$ and $J(a_1,\ldots,a_5)$ are 
related by the formula 
\begin{equation}\label{eq:I=cJ}
I(a_1,\ldots,a_5)=b\,J(a_1,\ldots,a_5)
\end{equation}
provided that $|p|$ is sufficiently small.  
To determine the constant $b$, we investigate the 
behavior of these two functions along the divisor $a_1a_2=1$. 
\par\medskip
We first consider the limit of $J(a_1,\ldots,a_5)$ as $a_1\to a_2^{-1}$. 
Since $\Gamma(a_1a_2;p,q)$ is written as 
\begin{equation*}
\Gamma(a_1a_2;p,q)=\frac{(pq/a_1a_2;p,q)_\infty}
{(1-a_1a_2)(pa_1a_2;p)_\infty(qa_1a_2;q)_\infty
(pqa_1a_2;p,q)_\infty},
\end{equation*}
we have 
\begin{equation}
\label{eq:lim(1-a1a2)Gamma(a1a2)}
\lim_{a_1\to a_2^{-1}}(1-a_1a_2)\Gamma(a_1a_2;p,q)
=\frac{1}{(p;p)_\infty(q;q)_\infty}.
\end{equation}
Using this, from \eqref{eq:defJ} we have 
\begin{equation}
\label{eq:calculation-J}
\begin{split}
\lim_{a_1\to a_2^{-1}}(1-a_1a_2)J(a)
&=
\frac{1}{(p;p)_\infty (q;q)_\infty}
\frac{\Gamma(a_2^{\pm2};p,q)}{\Gamma(a_2^{\pm1};p,q)}
\prod_{3\le i<j\le 5}\Gamma(a_ia_j;p,q)^2\ 
\Gamma(a_2^{\pm1}a_ia_j;p,q)\,
\\
&\quad\times
\Gamma(a_3a_4a_5;p,q)
\Gamma(a_2^{\pm1}a_3a_4a_5;p,q)
\prod_{k=3}^{5}\Gamma(a_k^2,a_2^{\pm1}a_k;p,q),
\end{split}
\end{equation}
where $a_5$ in the right-hand side should be understood as 
$a_{5}=\epsilon p^{1\over 2}q^{1\over 2}/a_1a_2a_3a_4$. 
Since we have 
$a_5\to \epsilon\,p^{1\over 2} q^{1\over 2}/a_3a_4$ in the limit $a_1\to a_2^{-1}$, 
from the definition \eqref{eq:def RGamma} of $\Gamma(u;p,q)$ 
we obtain 
\begin{equation*}
\begin{split}
\Gamma(a_3a_4a_5;p,q)&\to \Gamma(\epsilon p^{1\over 2}q^{1\over 2};p,q)=1,\\
\Gamma(a_2a_3a_4a_5,a_3a_4a_5/a_2;p,q)
&\to \Gamma(\epsilon p^{1\over 2}q^{1\over 2}a_2,\epsilon p^{1\over 2}q^{1\over 2}/a_2;p,q)=1
\end{split}
\end{equation*}
as $a_1\to a_2^{-1}$. 
Using this, 
\eqref{eq:calculation-J} implies 
\begin{lem} \label{lem:lim(1-a1a2)J(a)}
In the limit as $a_1\to a_2^{-1}$, the function $J(a)$ satisfies that 
\begin{equation*}
\begin{split}
\lim_{a_1\to a_2^{-1}}(1-a_1a_2)J(a)
&=
\frac{1}{(p;p)_\infty (q;q)_\infty}
\frac{\Gamma(a_2^{\pm2};p,q)}{\Gamma(a_2^{\pm1};p,q)}
\prod_{k=3}^{5}\Gamma(a_k^2,a_2^{\pm1}a_k;p,q)
\\
&\quad\times
\prod_{3\le i<j\le 5}\Gamma(a_ia_j;p,q)^2\ 
\Gamma(a_2^{\pm1}a_ia_j;p,q).
\end{split}
\end{equation*}
\end{lem}
\par
We next investigate the behavior of $I(a_1,\ldots,a_5)$ as $a_1\to a_2^{-1}$, 
assuming that $|p|$ is sufficiently small so that equality \eqref{eq:I=cJ} holds.  
Here, for convenience, we suppose that $|p|<|q|^{9}$ as in Theorem  \ref{thm:I=cJ}.  
We denote by $C_{r}(c)$ the positively oriented circle with center $c$ and radius $r$, i.e., 
$$
C_{r}(c)=\{z\in \mathbb{C}\,|\, |z-c|=r\}.
$$
Under the condition $a_{5}=\epsilon p^{1\over 2}q^{1\over 2}/a_1a_2a_3a_4$ with $|a_k|<1$ $(k=1,\ldots, 5)$, 
we consider $I(a)$ defined by \eqref{eq:I(a)} as the iterated integral
$$
I(a)=\frac{1}{(2\pi \sqrt{-1})^2}\int_{\mathbb{T}}
\Bigg(
\int_{\mathbb{T}}\Phi(a;z)\frac{dz_1}{z_1}\Bigg)
\frac{dz_2}{z_2},
$$
where the integrand $\Phi(a;z)$ defined by \eqref{eq:Phi} is written as 
\begin{equation}
\label{eq:Phi-2}
\Phi(a;z)=\frac{\prod_{k=1}^{5}
\Gamma(a_k z_2^{\pm1},a_k(z_1z_2)^{\pm1},
a_k(z_1z_2^2)^{\pm1};p,q)
}{
\Gamma(z_2^{\pm1},(z_1z_2)^{\pm1},
(z_1z_2^2)^{\pm1},
z_1^{\pm1},(z_1z_2^3)^{\pm1},
(z_1^2z_2^3)^{\pm1};p,q)
}.
\end{equation}
The function $\Phi(a;z)$ is rewritten as 
\begin{equation}
\label{eq:Phi-2+}
\Phi(a;z)=\Delta(z;p,q)\prod_{k=1}^{5}
\Gamma(a_k z_2^{\pm1},a_k(z_1z_2)^{\pm1},
a_k(z_1z_2^2)^{\pm1};p,q),
\end{equation}
where $\Delta(z;p,q)$ 
denotes the reciprocal of the denominator of $\Phi(a;z)$ in the expression \eqref{eq:Phi-2}. 
In fact, $\Delta(z;p,q)$ turns out to be a holomorphic function of  $z\in (\mathbb{C}^*)^2$, since 
it is computed as 
\begin{equation}
\label{eq:DeltaDelta}
\Delta(z;p,q)%
=\frac{1}{\Gamma(z_2^{\pm1},(z_1z_2)^{\pm1},
(z_1z_2^2)^{\pm1},z_1^{\pm1},(z_1z_2^3)^{\pm1},
(z_1^2z_2^3)^{\pm1};p,q)
}
=
\Delta(z;p)\Delta(z;q),
\end{equation}
where $\Delta(z;p)=z_1^{-3}z_2^{-5}
\theta(z_2,z_1z_2,z_1z_2^2,z_1,z_1z_2^3,z_1^2z_2^3;p)$ 
is just the elliptic Weyl denominator introduced by \eqref{eq:Delta(z;p)} in Section \ref{section:5}, 
and \eqref{eq:DeltaDelta} is immediately confirmed from \eqref{eq:RGamma2}. 
We also write $I(a)$ as 
\begin{equation}
\label{eq:I1(a,z2)}
I(a)=
\frac{1}{2\pi \sqrt{-1}}\int_{\mathbb{T}}
I_1(a;z_2)
\frac{dz_2}{z_2},
\quad\mbox{
where
}\quad
I_1(a;z_2)=\frac{1}{2\pi \sqrt{-1}}\int_{\mathbb{T}}\Phi(a;z)\frac{dz_1}{z_1}. 
\end{equation}
Since $\Phi(a;z)$ is rewritten from \eqref{eq:Phi-2} as 
\begin{equation}
\label{eq:Phi-3}
\Phi(a;z)=\frac{\prod_{k=1}^{5}\Gamma(a_kz_2^{\pm1};p,q)}
{\Gamma(z_2^{\pm1};p,q)}
\widetilde\Phi(a;z_1,z_2),
\end{equation}
where
\begin{equation*}
\widetilde\Phi(a;z_1,z_2)=
\frac{\prod_{k=1}^{5}
\Gamma(a_k(z_1z_2)^{\pm1},
a_k(z_1z_2^2)^{\pm1};p,q)
}{
\Gamma((z_1z_2)^{\pm1},
(z_1z_2^2)^{\pm1},
z_1^{\pm1},(z_1z_2^3)^{\pm1},
(z_1^2z_2^3)^{\pm1};p,q)
},
\end{equation*}
$I_1(a;z_2)$ is written as  
\begin{equation*}
I_1(a;z_2)=
\frac{\prod_{k=1}^{5}\Gamma(a_kz_2^{\pm1};p,q)}
{\Gamma(z_2^{\pm1};p,q)}
\frac{1}{2\pi \sqrt{-1}}\int_{\mathbb{T}}\widetilde\Phi(a;z_1,z_2)\frac{dz_1}{z_1}. 
\end{equation*}
\par
\begin{figure}[htbp]
\begin{center}
\includegraphics[width=405pt]{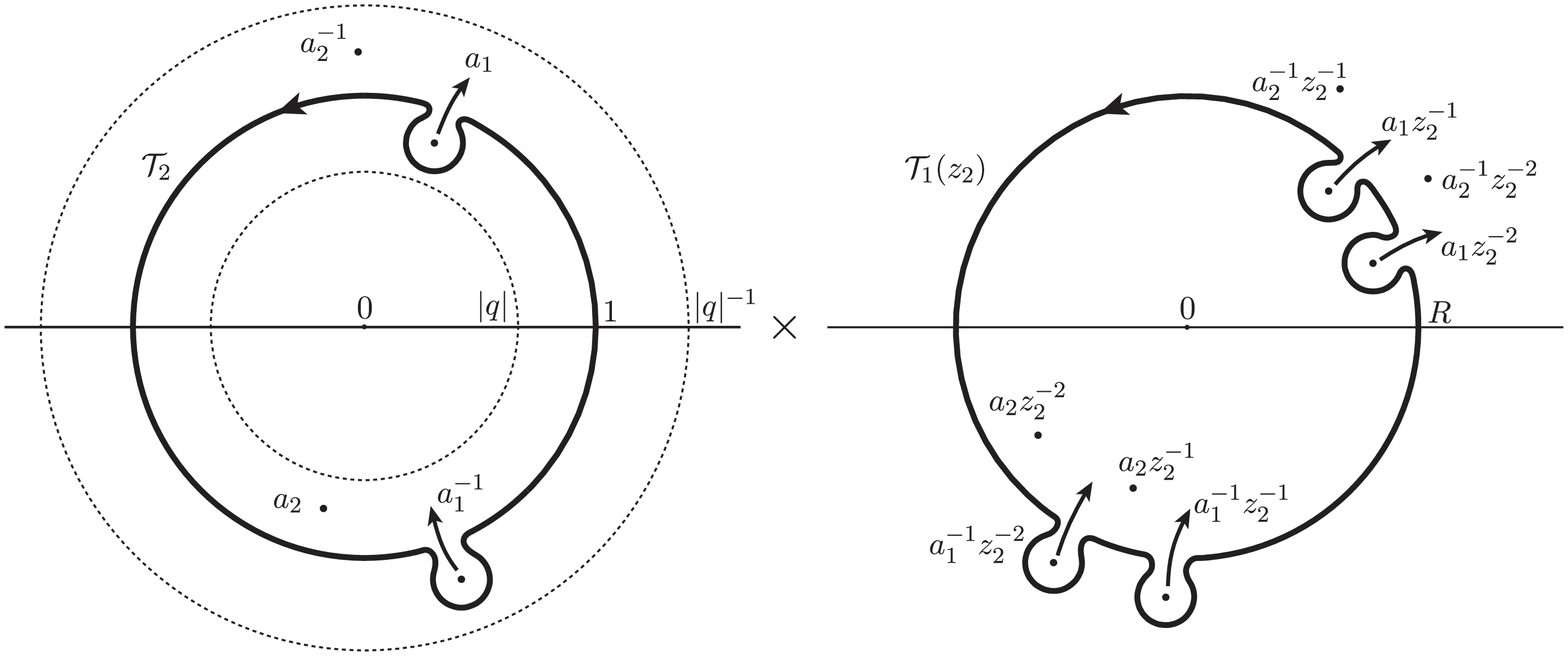}
\end{center}
\vspace{-15pt}
\begin{center}
\includegraphics[width=405pt]{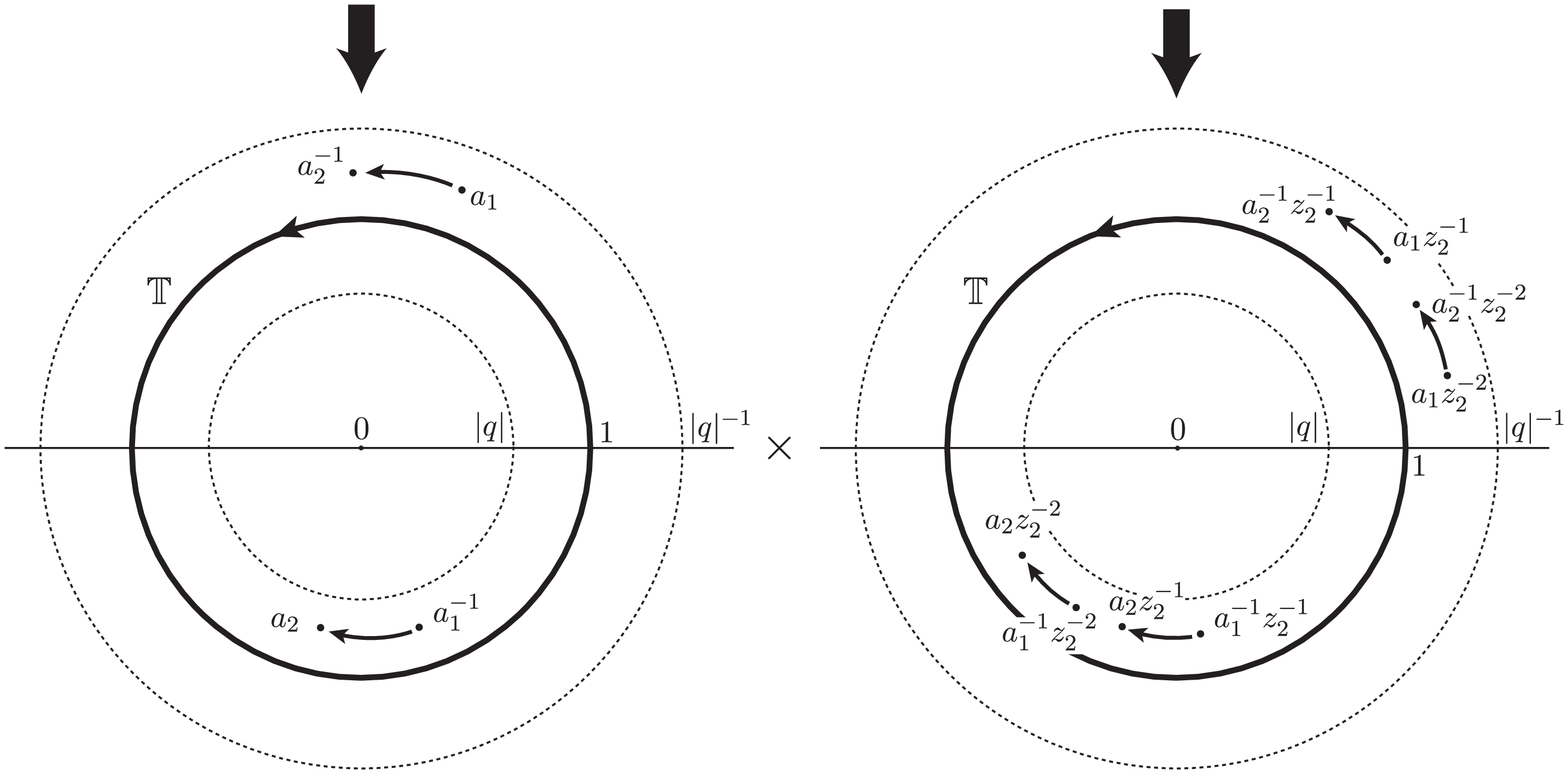}
\end{center}
\vspace{-17pt}
\caption{Deformation of the cycles}
\label{fig:complex planes}
\end{figure}

We now suppose $|a_i|< |a_1|<1$ $(i=2,\ldots,5)$. Then, the integral 
\begin{equation}
\label{eq:int tilde Phi}
\frac{1}{2\pi \sqrt{-1}}\int_{\mathbb{T}}\widetilde\Phi(a;z_1,z_2)\frac{dz_1}{z_1}
\end{equation}
defines a holomorphic function of $z_2$ on the domain 
$|a_1|^{1\over 2}<|z_2|<|a_1|^{-{1\over 2}}$. 
For any $z_2\in \mathbb{C}$ satisfying $|a_1|^2<|z_2|<|a_1|^{-2}$, 
it is easy to see that there exists $R>0$ such that the circle 
$C_R(0)$ in the $z_1$-plane keeps the points $z_1=a_1z_2^{-1}$, $a_1z_2^{-2}$ inside and 
$z_1=a_1^{-1}z_2^{-1}$, $a_1^{-1}z_2^{-2}$ outside. 
Then, the integral in \eqref{eq:int tilde Phi} is continued analytically to the domain 
$|a_1|^2<|z_2|<|a_1|^{-2}$
by deforming $\mathbb{T}=C_1(0)$ to $C_R(0)$. Hence 
\begin{equation*}
I_1(a;z_2)=
\frac{\prod_{k=1}^{5}\Gamma(a_kz_2^{\pm1};p,q)}
{\Gamma(z_2^{\pm1};p,q)}
\frac{1}{2\pi \sqrt{-1}}\int_{C_R(0)}\widetilde\Phi(a;z_1,z_2)\frac{dz_1}{z_1} 
\end{equation*}
defines a meromorphic function on the domain $|a_1|^2<|z_2|<|a_1|^{-2}$. 
\begin{lem}
\label{lem:I1(a,z2)}
Suppose $|a_i|< |a_1|<1$ $(i=2,\ldots,5)$.
The integral $I_1(a;z_2)$ is also expressed as 
\begin{equation}
\label{eq:I1(a,z2)-2}
\begin{split}
I_1(a;z_2)
&=\frac{1}{2\pi \sqrt{-1}}\int_{\mathcal{T}_1(z_2)}\Phi(a;z)\frac{dz_1}{z_1}
+
\frac{2\,\Gamma(a_1^2;p,q)
\prod_{k=2}^{5}\Gamma(a_ka_1^{\pm1},a_kz_2^{\pm1};p,q)
}{(p;p)_\infty(q;q)_\infty\Gamma(a_1^{\pm1},a_1^{-1}z_2^{\pm1};p,q)}
\\
&\qquad
\times
\Bigg(
\frac{
\prod_{k=2}^{5}
\Gamma(a_k(a_1z_2)^{\pm1};p,q)}
{\Gamma(z_2,
a_1^{-2}z_2^{-1},(a_1z_2^2)^{\pm1};p,q)}
+
\frac{\prod_{k=2}^{5}
\Gamma(a_k(a_1^{-1}z_2)^{\pm1};p,q)}
{\Gamma(z_2^{-1},
a_1^{-2}z_2,(a_1^{-1}z_2^2)^{\pm1};p,q)}
\Bigg),
\end{split}
\end{equation}
%
using  a cycle $\mathcal{T}_1(z_2)$ defined by the homology equivalence
\begin{equation}
\label{eq:C1(z2)}
\mathcal{T}_1(z_2)
\sim C_R(0)
-C_\varepsilon(a_1z_2^{-1})
+C_\varepsilon(a_1^{-1}z_2^{-1})
-C_\varepsilon(a_1z_2^{-2})
+C_\varepsilon(a_1^{-1}z_2^{-2}),
\end{equation}
where  
$\varepsilon>0$ is sufficiently small. 
\begin{figure}[htbp]
 \begin{center}
\includegraphics[height=155pt]{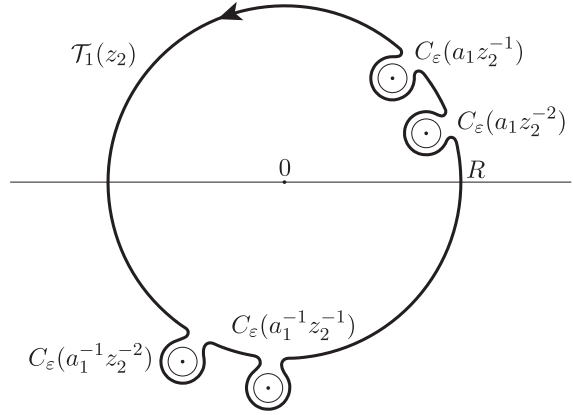}
\end{center}
\vspace{-20pt}
\caption{$z_1$-plane}
\label{fig:one}
\end{figure}
\vspace{-10pt}
\end{lem}
%
%
{\it Proof.} 
Since $C_R(0)\sim\mathcal{T}_1(z_2)
+C_\varepsilon(a_1z_2^{-1})
-C_\varepsilon(a_1^{-1}z_2^{-1})
+C_\varepsilon(a_1z_2^{-2})
-C_\varepsilon(a_1^{-1}z_2^{-2})$ by \eqref{eq:C1(z2)}, 
from 
\eqref{eq:I1(a,z2)}
we have 
\begin{align}
I_1(a;z_2)
&=\frac{1}{2\pi \sqrt{-1}}
\bigg(\int_{\mathcal{T}_1(z_2)}\Phi(a;z)\frac{dz_1}{z_1}
+\int_{C_\varepsilon(a_1z_2^{-1})}\Phi(a;z)\frac{dz_1}{z_1}-\int_{C_\varepsilon(a_1^{-1}z_2^{-1})}\Phi(a;z)\frac{dz_1}{z_1}
\nonumber\\[4pt]
&\hspace{85pt}
+\int_{C_\varepsilon(a_1z_2^{-2})}\Phi(a;z)\frac{dz_1}{z_1}-\int_{C_\varepsilon(a_1^{-1}z_2^{-2})}\Phi(a;z)\frac{dz_1}{z_1}
\bigg). 
\label{eq:I1(a,z2)-3}
\end{align}
Since $\Phi(a;z)$ is written as \eqref{eq:Phi-2+} with \eqref{eq:DeltaDelta}, 
using the formula \eqref{eq:lim(1-a1a2)Gamma(a1a2)}, we have 
%
\begin{align}
\label{eq:Res(z1=a1/z2)}
&\frac{1}{2\pi \sqrt{-1}}
\int_{C_\varepsilon(a_1z_2^{-1})}\Phi(a;z)\frac{dz_1}{z_1}=
\mathrm{Res}\Big(\Phi(a;z)\frac{dz_1}{z_1};
z_1=a_1z_2^{-1}\Big)
\nonumber\\
&\quad
=\lim_{z_1\to a_1z_2^{-1}}(z_1-a_1z_2^{-1})\frac{\Phi(a;z)}{z_1}
=\lim_{z_1\to a_1z_2^{-1}}(1-a_1z_1^{-1}z_2^{-1})\Phi(a;z)
\nonumber\\
&\quad
=\frac{\Gamma(a_1^2;p,q)}{(p;p)_\infty (q;q)_\infty}
\Delta(a_1z_2^{-1},z_2;p,q)
\prod_{k=2}^{5}
\Gamma(a_1^{\pm1}a_k;p,q)
\prod_{k=1}^{5}
\Gamma(a_kz_2^{\pm1}, a_k(a_1z_2)^{\pm1};p,q)
\nonumber\\
&\quad
=\frac{\Gamma(a_1^2;p,q)\prod_{k=2}^{5}
\Gamma(a_1^{\pm1}a_k;p,q)
\prod_{k=1}^{5}
\Gamma(a_kz_2^{\pm1}, a_k(a_1z_2)^{\pm1};p,q)}
{(p;p)_\infty (q;q)_\infty\Gamma(a_1^{\pm1},z_2^{\pm1},a_1^{\pm1}z_2^{\pm1},
(a_1z_2^2)^{\pm1},(a_1^2 z_2)^{\pm1};p,q)}
\nonumber\\
&\quad
=\frac{\Gamma(a_1^2;p,q)
\prod_{k=2}^{5}\Gamma(a_ka_1^{\pm1},a_kz_2^{\pm1},a_k(a_1z_2)^{\pm1};p,q)
}{(p;p)_\infty(q;q)_\infty\Gamma(a_1^{\pm1},a_1^{-1}z_2^{\pm1},z_2,
a_1^{-2}z_2^{-1},(a_1z_2^2)^{\pm1};p,q)}
\end{align}
and
%
\begin{align}
\label{eq:Res(z1=1/a1z2)}
&\frac{1}{2\pi \sqrt{-1}}
\int_{C_\varepsilon(a_1^{-1}z_2^{-1})}\Phi(a;z)\frac{dz_1}{z_1}
=\mathrm{Res}\Big(\Phi(a;z)\frac{dz_1}{z_1};
z_1=a_1^{-1}z_2^{-1}\Big)
\nonumber
\\&\quad
=\lim_{z_1\to a_1^{-1}z_2^{-1}}(z_1-a_1^{-1}z_2^{-1})\frac{\Phi(a;z)}{z_1}
=\lim_{z_1\to a_1^{-1}z_2^{-1}}\frac{-1}{a_1z_1z_2}(1-a_1z_1z_2)\Phi(a;z)
\nonumber
\\&\quad
=
-\frac{\Gamma(a_1^2;p,q)}{(p;p)_\infty (q;q)_\infty}
\Delta(a_1^{-1}z_2^{-1},z_2;p,q)
\prod_{k=2}^{5}
\Gamma(a_1^{\pm1}a_k;p,q)
\prod_{k=1}^{5}
\Gamma(a_kz_2^{\pm1}, a_k(a_1^{-1}z_2)^{\pm1};p,q)
\nonumber
\\&\quad
=
-\frac{\Gamma(a_1^2;p,q)\prod_{k=2}^{5}
\Gamma(a_1^{\pm1}a_k;p,q)
\prod_{k=1}^{5}
\Gamma(a_kz_2^{\pm1}, a_k(a_1^{-1}z_2)^{\pm1};p,q)}{(p;p)_\infty (q;q)_\infty
\Gamma(a_1^{\pm1},z_2^{\pm1},a_1^{\pm1}z_2^{\pm1},
(a_1^{-1}z_2^2)^{\pm1},(a_1^{-2}z_2)^{\pm1};p,q)}
\nonumber
\\&\quad
=
-\frac{\Gamma(a_1^2;p,q)
\prod_{k=2}^{5}\Gamma(a_ka_1^{\pm1},a_kz_2^{\pm1},a_k(a_1^{-1}z_2)^{\pm1};p,q)
}{(p;p)_\infty(q;q)_\infty\Gamma(a_1^{\pm1},a_1^{-1}z_2^{\pm1},z_2^{-1},
a_1^{-2}z_2,(a_1^{-1}z_2^2)^{\pm1};p,q)}.
\end{align}
Moreover, in the same way as above, we also have  
\begin{align}
\frac{1}{2\pi \sqrt{-1}}
\int_{C_\varepsilon(a_1z_2^{-2})}\Phi(a;z)\frac{dz_1}{z_1}
=
\frac{\Gamma(a_1^2;p,q)
\prod_{k=2}^{5}\Gamma(a_ka_1^{\pm1},a_kz_2^{\pm1},a_k(a_1^{-1}z_2)^{\pm1};p,q)
}{(p;p)_\infty(q;q)_\infty\Gamma(a_1^{\pm1},a_1^{-1}z_2^{\pm1},z_2^{-1},
a_1^{-2}z_2,(a_1^{-1}z_2^2)^{\pm1};p,q)},
\end{align}
which coincides with \eqref{eq:Res(z1=1/a1z2)} up to sign, and
\begin{align}
\frac{1}{2\pi \sqrt{-1}}
\int_{C_\varepsilon(a_1^{-1}z_2^{-2})}\Phi(a;z)\frac{dz_1}{z_1}
=-
\frac{\Gamma(a_1^2;p,q)
\prod_{k=2}^{5}\Gamma(a_ka_1^{\pm1},a_kz_2^{\pm1},a_k(a_1z_2)^{\pm1};p,q)
}{(p;p)_\infty(q;q)_\infty\Gamma(a_1^{\pm1},a_1^{-1}z_2^{\pm1},z_2,
a_1^{-2}z_2^{-1},(a_1z_2^2)^{\pm1};p,q)},
\label{eq:Res(z1=1/a1z2^2)}
\end{align}
which also coincides with \eqref{eq:Res(z1=a1/z2)} up to sign. 
Hence, applying \eqref{eq:Res(z1=a1/z2)}--\eqref{eq:Res(z1=1/a1z2^2)} to 
\eqref{eq:I1(a,z2)-3}, we see that $I_1(a;z_2)$ is expressed as \eqref{eq:I1(a,z2)-2} 
in Lemma \ref{lem:I1(a,z2)}. \qed
\begin{lem} 
\label{lem:I(a)2}
Suppose $|a_i|< |a_1|<1$ $(i=2,\ldots,5)$.
Then, it follows that 
\begin{equation}
\label{eq:I(a)2}
\begin{split}
I(a)
&=\frac{1}{(2\pi \sqrt{-1})^2}\int_{\mathcal{T}_2}
\bigg(\int_{\mathcal{T}_1(z_2)}\Phi(a;z)\frac{dz_1}{z_1}\bigg)
\frac{dz_2}{z_2}
\\&\quad
+
\frac{\Gamma(a_1^2;p,q)
\prod_{k=2}^{5}\Gamma(a_1^{\pm1}a_k;p,q)
}{(p;p)_\infty(q;q)_\infty\Gamma(a_1^{\pm1};p,q)}
\frac{1}{2\pi\sqrt{-1}}
\\&\qquad\times\Bigg(
2\int_{|z_2|=1}
\frac{
\prod_{k=2}^{5}
\Gamma(a_kz_2^{\pm1},a_k(a_1z_2)^{\pm1};p,q)}
{\Gamma(z_2,a_1^{-1}z_2^{\pm1},
a_1^{-2}z_2^{-1},(a_1z_2^2)^{\pm1};p,q)}
\frac{dz_2}{z_2}
\\&\quad\qquad
+
2\int_{|z_2|=1}
\frac{\prod_{k=2}^{5}
\Gamma(a_kz_2^{\pm1},a_k(a_1^{-1}z_2)^{\pm1};p,q)}
{\Gamma(z_2^{-1},a_1^{-1}z_2^{\pm1},
a_1^{-2}z_2,(a_1^{-1}z_2^2)^{\pm1};p,q)}
\frac{dz_2}{z_2}
\\&\quad\qquad
+
\int_{|z_1|=|a_1|^{-1}}
\frac{
\prod_{k=2}^{5}
\Gamma(a_k(a_1z_1)^{\pm1},a_k(a_1^2z_1)^{\pm1};p,q)
}{\Gamma(
z_1,
a_1z_1,a_1^{-2}z_1^{-1},
a_1^{-3}z_1^{-1},(a_1^3z_1^2)^{\pm1};p,q)}
\frac{dz_1}{z_1}
\\&\quad\qquad
+
\int_{|z_1|=|a_1|}
\frac{
\prod_{k=2}^{5}
\Gamma(a_k(a_1^{-1}z_1)^{\pm1},a_k(a_1^{-2}z_1)^{\pm1};p,q)
}{\Gamma(
z_1^{-1},
a_1z_1^{-1},a_1^{-2}z_1,
a_1^{-3}z_1,(a_1^{-3}z_1^2)^{\pm1};p,q)}
\frac{dz_1}{z_1}
\Bigg),
\end{split}
\end{equation}
using a cycle $\mathcal{T}_2$ defined by the homological equivalence 
\begin{equation}
\label{eq:calC2}
\mathcal{T}_2\sim\mathbb{T}
-C_\varepsilon(a_1)
+C_\varepsilon(a_1^{-1}),
\end{equation}
where $\varepsilon>0$ is sufficiently small. 
\begin{figure}[htbp]
\begin{center}
\includegraphics[height=155pt]{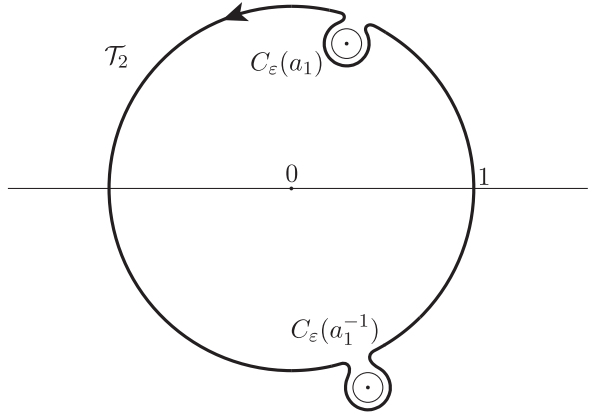}
\end{center}
\vspace{-20pt}
\caption{$z_2$-plane}
\end{figure}
\vspace{-15pt}
\end{lem}
{\it Proof.} From Lemma \ref{lem:I1(a,z2)}, \eqref{eq:I1(a,z2)} implies that
\begin{equation}
\label{eq:I(a)3}
\begin{split}
I(a)&=
\frac{1}{2\pi \sqrt{-1}}\int_{\mathbb{T}}
I_1(a;z_2)
\frac{dz_2}{z_2}\\
&=
\frac{1}{(2\pi \sqrt{-1})^2}\int_{\mathbb{T}}\bigg(
\int_{\mathcal{T}_1(z_2)}\Phi(a;z)\frac{dz_1}{z_1}
\bigg)\frac{dz_2}{z_2}
+
\frac{\Gamma(a_1^2;p,q)
\prod_{k=2}^{5}\Gamma(a_ka_1^{\pm1};p,q)
}{(p;p)_\infty(q;q)_\infty\Gamma(a_1^{\pm1};p,q)}\frac{1}{2\pi \sqrt{-1}}
\\
&\hspace{150pt}\times
\Bigg(2\int_{\mathbb{T}}
\frac{
\prod_{k=2}^{5}
\Gamma(a_kz_2^{\pm1},a_k(a_1z_2)^{\pm1};p,q)}
{\Gamma(a_1^{-1}z_2^{\pm1},z_2,
a_1^{-2}z_2^{-1},(a_1z_2^2)^{\pm1};p,q)}
\frac{dz_2}{z_2}\\
&\hspace{150pt}\quad+
2\int_{\mathbb{T}}
\frac{\prod_{k=2}^{5}
\Gamma(a_kz_2^{\pm1},a_k(a_1^{-1}z_2)^{\pm1};p,q)}
{\Gamma(a_1^{-1}z_2^{\pm1},z_2^{-1},
a_1^{-2}z_2,(a_1^{-1}z_2^2)^{\pm1};p,q)}
\frac{dz_2}{z_2}
\Bigg).
\end{split}
\end{equation}
Since 
$\mathbb{T}\sim \mathcal{T}_2+C_\varepsilon(a_1)-C_\varepsilon(a_1^{-1})$ 
by \eqref{eq:calC2}, 
the initial term of the right-hand side of \eqref{eq:I(a)3} is 
\begin{equation}
\label{eq:initial term SS}
\begin{split}
&
\frac{1}{(2\pi \sqrt{-1})^2}
\int_{\mathbb{T}}\bigg(
\int_{\mathcal{T}_1(z_2)}\Phi(a;z)\frac{dz_1}{z_1}
\bigg)\frac{dz_2}{z_2}
=
\frac{1}{(2\pi \sqrt{-1})^2}\Bigg[
\int_{\mathcal{T}_2}\bigg(
\int_{\mathcal{T}_1(z_2)}\Phi(a;z)\frac{dz_1}{z_1}
\bigg)\frac{dz_2}{z_2} \\[2pt]
&\hspace{70pt}
+
\int_{C_\varepsilon(a_1)}\bigg(
\int_{\mathcal{T}_1(z_2)}\Phi(a;z)\frac{dz_1}{z_1}
\bigg)\frac{dz_2}{z_2}
-
\int_{C_\varepsilon(a_1^{-1})}\bigg(
\int_{\mathcal{T}_1(z_2)}\Phi(a;z)\frac{dz_1}{z_1}
\bigg)\frac{dz_2}{z_2}
\Bigg].
\end{split}
\end{equation}
Since $\Phi(a;z)$ is written as \eqref{eq:Phi-3},  
using \eqref{eq:lim(1-a1a2)Gamma(a1a2)},
the second term in the right-hand side of \eqref{eq:initial term SS} is calculated as 
\begin{equation}
\label{eq:second term SS}
\begin{split}
&\frac{1}{(2\pi \sqrt{-1})^2}
\int_{C_\varepsilon(a_1)}\bigg(
\int_{\mathcal{T}_1(z_2)}\Phi(a;z)\frac{dz_1}{z_1}
\bigg)\frac{dz_2}{z_2}\\
&\quad=
\frac{1}{(2\pi \sqrt{-1})^2}
\int_{C_\varepsilon(a_1)}
\frac{\prod_{k=1}^{5}\Gamma(a_kz_2^{\pm1};p,q)}
{\Gamma(z_2^{\pm1};p,q)}
\bigg(
\int_{\mathcal{T}_1(z_2)}
\widetilde\Phi (a;z_1,z_2)
\frac{dz_1}{z_1}\bigg)\frac{dz_2}{z_2}
\\
&\quad=\frac{1}{2\pi \sqrt{-1}}\,
\mathrm{Res}\Bigg(\frac{\prod_{k=1}^{5}\Gamma(a_kz_2^{\pm1};p,q)}
{\Gamma(z_2^{\pm1};p,q)}
\bigg(
\int_{\mathcal{T}_1(z_2)}
\widetilde\Phi (a;z_1,z_2)
\frac{dz_1}{z_1}\bigg)\frac{dz_2}{z_2};
z_2=a_1\Bigg)
\\
&\quad=
\mathrm{Res}\bigg(\frac{\prod_{k=1}^{5}\Gamma(a_kz_2^{\pm1};p,q)}
{\Gamma(z_2^{\pm1};p,q)}
\frac{dz_2}{z_2};
z_2=a_1\bigg)
\frac{1}{2\pi \sqrt{-1}}
\int_{\mathcal{T}_1(a_1)}
\widetilde\Phi (a;z_1,a_1)
\frac{dz_1}{z_1}\\
&\quad=
\bigg(\lim_{z_2\to a_1}\frac{z_2-a_1}{z_2}\frac{\prod_{k=1}^{5}\Gamma(a_kz_2^{\pm1};p,q)}
{\Gamma(z_2^{\pm1};p,q)}
\bigg)
\frac{1}{2\pi \sqrt{-1}}
\int_{\mathcal{T}_1(a_1)}
\widetilde\Phi (a;z_1,a_1)
\frac{dz_1}{z_1}\\
&\quad=
\frac{\Gamma(a_1^2;p,q)\prod_{k=2}^{5}\Gamma(a_ka_1^{\pm1};p,q)
}{(p;p)_\infty(q;q)_\infty\Gamma(a_1^{\pm1};p,q)}
\frac{1}{2\pi \sqrt{-1}}
\int_{\mathcal{T}_1(a_1)}
\widetilde\Phi (a;z_1,a_1)
\frac{dz_1}{z_1}.
\end{split}
\end{equation}
Since the function 
$$
\widetilde\Phi (a;z_1,a_1)=\frac{
\prod_{k=2}^{5}
\Gamma(a_k(a_1z_1)^{\pm1},a_k(a_1^2z_1)^{\pm1};p,q)
}{\Gamma(
z_1,
a_1z_1,a_1^{-2}z_1^{-1},
a_1^{-3}z_1^{-1},(a_1^3z_1^2)^{\pm1};p,q)}
$$
of $z_1$ is now holomorphic at the points $z_1=1,a_1^{-1},a_1^{-2},a_1^{-3}$, 
which are avoided by the contour $\mathcal{T}_1(a_1)$ of the integral, 
we can deform the contour $\mathcal{T}_1(a_1)$ to 
the circle $\{z_1\in \mathbb{C}\,|\, 
|z_1|=|a_1|^{-1}\}
$ across these points, provided 
$|a_i|< |a_1|<1$ $(i=2,\ldots,5)$.
This means that 
\begin{equation}
\label{eq:second term SS-2}
\int_{\mathcal{T}_1(a_1)}
\widetilde\Phi (a;z_1,a_1)
\frac{dz_1}{z_1}
=
\int_{|z_1|=|a_1|^{-1}}
\frac{
\prod_{k=2}^{5}
\Gamma(a_k(a_1z_1)^{\pm1},a_k(a_1^2z_1)^{\pm1};p,q)
}{\Gamma(
z_1,
a_1z_1,a_1^{-2}z_1^{-1},
a_1^{-3}z_1^{-1},(a_1^3z_1^2)^{\pm1};p,q)}
\frac{dz_1}{z_1}.
\end{equation}
From \eqref{eq:second term SS} and \eqref{eq:second term SS-2}, we therefore obtain 
\begin{equation}
\label{eq:second term SS-3}
\begin{split}
\frac{1}{(2\pi \sqrt{-1})^2}
&\int_{C_\varepsilon(a_1)}\bigg(
\int_{\mathcal{T}_1(z_2)}\Phi(a;z)\frac{dz_1}{z_1}
\bigg)\frac{dz_2}{z_2}
=
\frac{\Gamma(a_1^2;p,q)\prod_{k=2}^{5}\Gamma(a_ka_1^{\pm1};p,q)
}{(p;p)_\infty(q;q)_\infty\Gamma(a_1^{\pm1};p,q)}\\[2pt]
&\hspace{50pt}
\times
\frac{1}{2\pi \sqrt{-1}}
\int_{|z_1|=|a_1|^{-1}}
\frac{
\prod_{k=2}^{5}
\Gamma(a_k(a_1z_1)^{\pm1},a_k(a_1^2z_1)^{\pm1};p,q)
}{\Gamma(
z_1,
a_1z_1,a_1^{-2}z_1^{-1},
a_1^{-3}z_1^{-1},(a_1^3z_1^2)^{\pm1};p,q)}
\frac{dz_1}{z_1}.
\end{split}
\end{equation}
In the same way as above, 
the third term in the right-hand side of \eqref{eq:initial term SS} is also calculated as 
\begin{equation}
\label{eq:third term SS}
\begin{split}
\frac{1}{(2\pi \sqrt{-1})^2}
&\int_{C_\varepsilon(a_1^{-1})}\bigg(
\int_{\mathcal{T}_1(z_2)}\Phi(a;z)\frac{dz_1}{z_1}
\bigg)\frac{dz_2}{z_2}
=
-\frac{\Gamma(a_1^2;p,q)\prod_{k=2}^{5}\Gamma(a_ka_1^{\pm1};p,q)
}{(p;p)_\infty(q;q)_\infty\Gamma(a_1^{\pm1};p,q)}\\[2pt]
&\hspace{50pt}
\times
\frac{1}{2\pi \sqrt{-1}}
\int_{|z_1|=|a_1|}
\frac{
\prod_{k=2}^{5}
\Gamma(a_k(a_1^{-1}z_1)^{\pm1},a_k(a_1^{-2}z_1)^{\pm1};p,q)
}{\Gamma(
z_1^{-1},
a_1z_1^{-1},a_1^{-2}z_1,
a_1^{-3}z_1,(a_1^{-3}z_1^2)^{\pm1};p,q)}
\frac{dz_1}{z_1}.
\end{split}
\end{equation}
Hence, combining \eqref{eq:I(a)3}, \eqref{eq:initial term SS}, \eqref{eq:second term SS-3} and \eqref{eq:third term SS}, we obtain the expression \eqref{eq:I(a)2} of $I(a)$ in Lemma \ref{lem:I(a)2}.
\qed
\begin{lem}
\label{lem:lim(1-a1a2)I(a)}
Suppose that 
$|a_i|< |a_2|<1$ $(i=3,\ldots,5)$
and $|q|<|a_2|<1$. Then it follows that  
\begin{equation}
\label{eq:lim(1-a1a2)I(a)}
\begin{split}
\lim_{a_1\to a_2^{-1}}(1-a_1a_2)I(a)
&=
\frac{12}{(p;p)_\infty^3(q;q)_\infty^3}
\frac{\Gamma(a_2^{\pm2};p,q)}
{\Gamma(a_2^{\pm1};p,q)}
\prod_{k=3}^{5}
\Gamma(a_k^2,a_2^{\pm1}a_k;p,q)
\\
&\qquad\times\prod_{3\le i<j\le 5}\Gamma(a_ia_j;p,q)^2\ 
\Gamma(a_2^{\pm1}a_ia_j;p,q).
\end{split}
\end{equation}
\end{lem}
{\it Proof.} 
Before taking the limit $a_1\to a_2^{-1}$ for $I(a)$, 
we need to extend $I(a)$ by analytic continuation to the function of $a_1$ on the domain $1\le |a_1|<|q|^{-1}$, 
provided $1<|a_2^{-1}|<|q|^{-1}$. 
For $a_1$ satisfying $1< |a_1|<|q|^{-1}$, the integral of the first term of \eqref{eq:I(a)2} in Lemma \ref{lem:I(a)2} can be extended by 
\begin{equation}
\label{eq:first term SS}
\frac{1}{(2\pi \sqrt{-1})^2}\int_{\mathcal{T}_2}
\Bigg(\int_{\mathcal{T}_1(z_2)}\Phi(a;z)\frac{dz_1}{z_1}\Bigg)
\frac{dz_2}{z_2}
=
\frac{1}{(2\pi \sqrt{-1})^2}\int_{\mathbb{T}}
\Bigg(\int_{\mathbb{T}}\Phi(a;z)\frac{dz_1}{z_1}\Bigg)
\frac{dz_2}{z_2}
\end{equation}
deforming the cycle $\mathcal{T}_2\times \mathcal{T}_1(z_2)$ to $\mathbb{T}\times \mathbb{T}$ of the integral 
as 
Figure \ref{fig:complex planes}.\\[-13pt]
%

%
\noindent
Since the integral \eqref{eq:first term SS} as a function of $a_1$ 
is regular at $a_1=a_2^{-1}$, 
the integral \eqref{eq:first term SS} itself has a finite limit as $a_1\to a_2^{-1}$.  
Hence, using \eqref{eq:lim(1-a1a2)Gamma(a1a2)},  \eqref{eq:I(a)2} implies that
\begin{equation}
\label{eq:lim(1-a1a2)I(a)2}
\begin{split}
\lim_{a_1\to a_2^{-1}}(1-a_1a_2)I(a)
&=
\frac{\Gamma(a_2^{\pm 2};p,q)
\prod_{k=3}^{5}\Gamma(a_2^{\pm1}a_k;p,q)
}{(p;p)_\infty^2(q;q)_\infty^2\Gamma(a_2^{\pm1};p,q)}
\frac{1}{2\pi\sqrt{-1}}
\\
%
%
%
%
&\qquad\mbox{}\times\Bigg(
2\int_{|z_2|=1}
\frac{
\prod_{k=3}^{5}
\Gamma(a_kz_2^{\pm1},a_k(a_2^{-1}z_2)^{\pm1};p,q)}
{\Gamma((a_2^{-1}z_2^2)^{\pm1};p,q)}
\frac{dz_2}{z_2}
\\
&\quad\qquad
+
2\int_{|z_2|=1}
\frac{
\prod_{k=3}^{5}
\Gamma(a_kz_2^{\pm1},a_k(a_2z_2)^{\pm1};p,q)}
{\Gamma((a_2z_2^2)^{\pm1};p,q)}
\frac{dz_2}{z_2}
\\
&\quad\qquad
+
\int_{|z_1|=|a_2|}
\frac{
\prod_{k=3}^{5}
\Gamma(a_k(a_2^{-1}z_1)^{\pm1},a_k(a_2^{-2}z_1)^{\pm1};p,q)
}{\Gamma(a_2^{-3}z_1^2)^{\pm1};p,q)}
\frac{dz_1}{z_1}
\\
&\quad\qquad
+
\int_{|z_1|=|a_2|^{-1}}
\frac{
\prod_{k=3}^{5}
\Gamma(a_k(a_2z_1)^{\pm1},a_k(a_2^2z_1)^{\pm1};p,q)
}{\Gamma((a_2^{3}z_1^2)^{\pm1};p,q)}
\frac{dz_1}{z_1}
\Bigg),
\end{split}
\end{equation}
where $a_5$ in the right-hand side should be understood as 
$a_{5}=\epsilon p^{1\over 2}q^{1\over 2}/a_1a_2a_3a_4
\to \epsilon p^{1\over 2}q^{1\over 2}/a_3 a_4$. 
Applying the variable changes 
\begin{equation*}
a_2^{-{1\over 2}}z_2=w,
\quad 
a_2^{{1\over 2}}z_2=w,
\quad 
a_2^{-\frac{3}{2}}z_1=w,
\quad 
a_2^{\frac{3}{2}}z_1=w
\end{equation*}
to the integrals in \eqref{eq:lim(1-a1a2)I(a)2}, respectively from the top, we have 
\begin{equation}
\label{eq:lim(1-a1a2)I(a)3}
\begin{split}
&\lim_{a_1\to a_2^{-1}}(1-a_1a_2)I(a)
=
\frac{\Gamma(a_2^{\pm 2};p,q)
\prod_{k=3}^{5}\Gamma(a_2^{\pm1}a_k;p,q)
}{(p;p)_\infty^2(q;q)_\infty^2\Gamma(a_2^{\pm1};p,q)}
\frac{1}{2\pi\sqrt{-1}}\\
&\hspace{120pt}\times
\Bigg(
3\int_{|w|=|a_2|^{1\over 2}}
\frac{
\prod_{k=3}^{5}
\Gamma(a_ka_2^{1\over 2}w^{\pm1},a_ka_2^{-{1\over 2}}w^{\pm1};p,q)}
{\Gamma(w^{\pm2};p,q)}
\frac{dw}{w}
\\
&\hspace{140pt}
+
3\int_{|w|=|a_2|^{-{1\over 2}}}
\frac{
\prod_{k=3}^{5}
\Gamma(a_ka_2^{1\over 2}w^{\pm1},a_ka_2^{-{1\over 2}}w^{\pm1};p,q)}
{\Gamma(w^{\pm2};p,q)}
\frac{dw}{w}
\Bigg)\\
&=\frac{\Gamma(a_2^{\pm 2};p,q)
\prod_{k=3}^{5}\Gamma(a_2^{\pm1}a_k;p,q)
}{(p;p)_\infty^2(q;q)_\infty^2\Gamma(a_2^{\pm1};p,q)}
\frac{6}{2\pi\sqrt{-1}}
\int_{|w|=1}
\frac{
\prod_{k=3}^{5}
\Gamma(a_ka_2^{1\over 2}w^{\pm1},a_ka_2^{-{1\over 2}}w^{\pm1};p,q)}
{\Gamma(w^{\pm2};p,q)}
\frac{dw}{w},
\end{split}
\end{equation}
whose contour $|w|=1$ is deformed from the contours $|w|=|a_2|^{{1\over 2}}$ or $|w|=|a_2|^{-{1\over 2}}$,
since its integrand has no poles in the annulus $|a_2|^{{1\over 2}}\le |w|\le |a_2|^{-{1\over 2}}$, 
provided 
$|a_i|< |a_2|<1$ $(i=3,\ldots,5)$. 
Here the integral 
$$
\frac{1}{2\pi\sqrt{-1}}
\int_{|w|=1}
\frac{
\prod_{k=3}^{5}
\Gamma(a_ka_2^{1\over 2}w^{\pm1},a_ka_2^{-{1\over 2}}w^{\pm1};p,q)}
{\Gamma(w^{\pm2};p,q)}
\frac{dw}{w}
$$
coincides with the elliptic integral \eqref{eq:eAW} of type $BC_1$
for specific parameters
\begin{equation*}
a_2^{{1\over 2}} a_3, \ \ 
a_2^{-{1\over 2}} a_3, \ \ 
a_2^{{1\over 2}} a_4, \ \ 
a_2^{-{1\over 2}} a_4, \ \ 
a_2^{{1\over 2}} a_5, \ \ 
a_2^{-{1\over 2}} a_5 \ \ 
\end{equation*}
satisfying the balancing condition 
$\prod_{k=3}^5(a_2^{1\over 2}a_k)(a_2^{-{1\over 2}}a_k)=(a_3a_4a_5)^2=pq$, and is evaluated as 
\begin{equation*}
\frac{2\,
\Gamma(a_3^2,a_4^2,a_5^2,a_2^{\pm1}a_3a_4,a_2^{\pm1}a_3a_5,a_2^{\pm1}a_4a_5;p,q)
\Gamma(a_3a_4,a_3a_5,a_4a_5;p,q)^2
}{(p;p)_\infty(q;q)_\infty},
\end{equation*}
which is confirmed from the right-hand side of \eqref{eq:eAW}. 
Therefore, applying this to \eqref{eq:lim(1-a1a2)I(a)3}, 
we eventually obtain \eqref{eq:lim(1-a1a2)I(a)} in Lemma \ref{lem:lim(1-a1a2)I(a)}.  
This completes the proof. \qed
\medskip
Lastly we determine the constant $b$ independent of $a$ appearing in \eqref{eq:I=cJ}. 
Comparing both sides of \eqref{eq:I=cJ} as $a_1\to a_2^{-1}$ using Lemmas \ref{lem:lim(1-a1a2)J(a)} and \ref{lem:lim(1-a1a2)I(a)}, 
we obtain 
$$
b=\frac{I(a)}{J(a)}=\lim_{a_1\to a_2^{-1}}\frac{(1-a_1a_2)I(a)}{(1-a_1a_2)J(a)}
=\frac{12}{(p;p)_\infty^2(q;q)_\infty^2}.
$$
\section*{Acknowledgements} 
The authors are grateful to the anonymous referee for valuable comments that have helped them
improve the manuscript. This work is supported by JSPS Kakenhi Grants  
(B)15H03626 and (C)18K03339. 


{\footnotesize
\begin{thebibliography}{99}
%
%
\bibitem{AW}
 R.~Askey and J. Wilson: Some basic hypergeometric orthogonal polynomials that generalize Jacobi polynomials, 
 Mem.~Amer.~Math.~Soc.~{\bf 54} (1985), no.~319, iv+55 pp. 
%
\bibitem{Gu1994}
R.~A.~Gustafson: Some $q$-beta and Mellin--Barnes integrals on compact Lie groups and Lie algebras, 
Trans. Amer.~Math.~Soc.~{\bf 341} (1994), 69--119.
%
\bibitem{Gu1990}
R.~A.~Gustafson: A summation theorem for hypergeometric series very-well-poised on $G_2$, 
SIAM J.~Math. Anal.~{\bf 21} (1990), 510--522.
%
%
\bibitem{Ito}
M.~Ito: Askey--Wilson type integrals associated with root systems, Ramanujan J.~{\bf 12} (2006), 131--151.
%
\bibitem{IN2019}
M.~Ito and M.~Noumi: 
A determinant formula associated with the elliptic hypergeometric integrals of type $BC_n$, 
J.~Math.~Phys.~{\bf 60} (2019), 071705, 31 pp.
%
\bibitem{IN2018}
M.~Ito and M.~Noumi: 
Connection formula for the Jackson integral of type $A_n$ and elliptic Lagrange interpolation, 
SIGMA Symmetry Integrability Geom. Methods Appl.~{\bf 14} (2018), Paper No.~077, 42 pp.
%
\bibitem{INSumBC}
M.~Ito and M.~Noumi:
Derivation of a $BC_n$ elliptic summation formula via the fundamental invariants, 
Constr.~Approx. {\bf 45} (2017), 33--46.
\bibitem{INIntBC}
M.~Ito and M.~Noumi: Evaluation of the $BC_n$ elliptic Selberg integral 
via the fundamental invariants,
Proc.~Amer.~Math.~Soc.~{\bf 145} (2017), 689--703. 
%
\bibitem{INSlaterBC}
M.~Ito and M.~Noumi: A generalization of the Sears--Slater transformation 
and elliptic  Lagrange interpolation of type $BC_n$, 
Adv.~in Math.~{\bf 229} (2016), 361--380.
%
\bibitem{Mac}
I.~G.~Macdonald: Affine Hecke algebras and orthogonal polynomials, 
Cambridge Tracts in Mathematics, 157. Cambridge University Press, Cambridge, 2003. x+175 pp.
%
\bibitem{Ra2010}
E.~M.~Rains: Transformations of elliptic hypergeometric integrals. Ann.~of Math.~(2) {\bf 171} (2010), 169--243.
%
\bibitem{SV2011}
V.~P.~Spiridonov and G.~S.~Vartanov: Elliptic hypergeometry of supersymmetric dualities, Comm.~Math.~Phys. {\bf 304} (2011), 797--874.

\bibitem{SV2010}
V.~P.~Spiridonov and G.~S.~Vartanov: 
Superconformal indices for 
$\mathcal{N}=1$ theories with multiple duals, Nuclear Phys.~B {\bf 824} (2010), 192--216. [arXiv:0811.1909]
%
\bibitem{S}
V.~P.~Spiridonov: Short proofs of the elliptic beta integrals. Ramanujan J.~{\bf 13} (2007), 265--283. 
%
\bibitem{vDS}
J.~F.~van Diejen and V.~P.~Spiridonov: Elliptic Selberg integrals. Internat. Math. Res. Notices 2001. 1083--1110.


\end{thebibliography}
}
\end{document}